\documentclass{article}

\usepackage{arxiv}

\usepackage{todonotes}
\usepackage{hhline}

\usepackage{tabularx,tabu}
\usepackage[ruled,linesnumbered]{algorithm2e}
\usepackage{etoolbox}
\usepackage{algpseudocode}
\usepackage{etoolbox}
\usepackage{multirow, makecell}
\usepackage{framed}
\usepackage[export]{adjustbox}
\usepackage{mathtools}
\usepackage{enumitem}
\usepackage{amssymb}
\usepackage{amsthm}
\usepackage{amsmath}
\usepackage{upgreek}
\usepackage{amssymb}
\usepackage{commath}
\usepackage{comment}
\usepackage[utf8]{inputenc} 
\usepackage[T1]{fontenc}    
\usepackage{url}            
\usepackage{booktabs}       
\usepackage{amsfonts}       
\usepackage{nicefrac}       
\usepackage{microtype}      
\usepackage{lipsum}		
\usepackage{graphicx}
\usepackage{subcaption}
\usepackage{float}
\usepackage{lscape}
\usepackage{natbib}
\usepackage{doi}
\usepackage{color,soul}
\usepackage{titlesec}
\titlespacing\subsubsection{0pt}{12pt plus 4pt minus 2pt}{0pt plus 2pt minus 2pt}

\newtheorem{example}{Example}
\newtheorem{definition}{Definition}

\usepackage{caption}
\SetKwInOut{KwInput}{Input}
\SetKwInOut{KwOutput}{Output}
\SetKwInOut{KwInit}{Initialization}

\def\R{\mathbb{R}}

\def\bc{\boldsymbol{c}}
\def\a{\boldsymbol{a}}

\def\d{\boldsymbol{d}}
\def\h{\boldsymbol{h}}

\def\T{\boldsymbol{T}}
\def\W{\boldsymbol{W}}

\def\z{\boldsymbol{z}}
\def\Z{\mathcal{Z}}
\def\barZ{\bar{\Z}}

\def\x{\boldsymbol{x}}
\def\X{\mathcal{X}}

\def\y{\boldsymbol{y}}
\def\Y{\mathcal{Y}}
\def\K{\mathcal{K}}
\def\L{\mathcal{L}}

\def\N{\mathcal{N}}
\def\w{\boldsymbol{w}}
\def\1{\boldsymbol{1}}
\def\0{\boldsymbol{0}}

\def\s{\boldsymbol{s}}
\def\a{\boldsymbol{a}}

\def\P{\mathcal{P}}
\def\p{\boldsymbol{p}}
\def\q{\boldsymbol{q}}
\def\hq{\Hat{q}}

\title{Machine Learning for $K$-adaptability \\ in Two-Stage Robust Optimization}
\author{Esther Julien\\ 
        TU Delft\\
        \texttt{e.a.t.julien@tudelft.nl} \\
        \And
        Krzysztof Postek\\
        Independent \\ Researcher\\
        \texttt{}\\
        \And
        \c{S}. \.{I}lker Birbil\\
        University of Amsterdam\\
        \texttt{s.i.birbil@uva.nl}}

\renewcommand{\headeright}{}
\renewcommand{\undertitle}{}
\renewcommand{\shorttitle}{Machine Learning for $K$-adaptability}

\hypersetup{
pdftitle={Machine Learning for $K$-adaptability in Two-stage Robust Optimization},
pdfsubject={math.OC},
pdfauthor={Esther Julien, Krzysztof Postek, \c{S}. \.{I}lker Birbil},
pdfkeywords={node selection, clustering, two-stage robust optimization, K-adaptability, machine learning, tree search}}

\begin{document}
\maketitle

\begin{abstract}
    Two-stage robust optimization problems constitute one of the hardest optimization problem classes. One of the solution approaches to this class of problems is $K$-adaptability. This approach simultaneously seeks the best partitioning of the uncertainty set of scenarios into $K$ subsets, and optimizes decisions corresponding to each of these subsets. In general case, it is solved using the $K$-adaptability branch-and-bound algorithm, which requires exploration of exponentially-growing solution trees. To accelerate finding high-quality solutions in such trees, we propose a machine learning-based node selection strategy. In particular, we construct a feature engineering scheme based on general two-stage robust optimization insights that allows us to train our machine learning tool on a database of resolved B\&B trees, and to apply it as-is to problems of different sizes and/or types. We experimentally show that using our learned node selection strategy outperforms a vanilla, random node selection strategy when tested on problems of the same type as the training problems, also in case the $K$-value or the problem size differs from the training ones.
\end{abstract}


\section{Introduction} \label{sec:introduction}
Many optimization problems are affected by data uncertainty caused by errors in the forecast, implementation, or measurement. Robust optimization (RO) is one of the key paradigms to solve such problems, where the goal is to find an optimal solution among the ones that remain feasible for all data realizations within an \emph{uncertainty set} \citep{ben2009robust}. This set includes all \emph{reasonable} data outcomes.

A specific class of RO problems comprises two-stage robust optimization (2SRO) problems in which some decisions are implemented \emph{before} the uncertain data is known (here-and-now decisions), and other decisions are implemented \emph{after} the data is revealed (wait-and-see decisions). Such a problem can be formulated as
\begin{equation}
\min_{\x \in \X} \; \max_{\z \in \Z} \; \min_{\y \in \Y} \big \{\bc(\z)^\intercal \x + \d(\z)^\intercal \y: \T(\z)\x + \W(\z)\y \leq \h(\z), \; \forall \z \in \Z \big \}, \label{eq:sec1:two.stage}
\end{equation}
where $\x \in \X \subseteq \R^{N_x}$ and $\y \in \Y \subseteq \R^{N_y}$ are the here-and-now and wait-and-see decisions, respectively, and $\z$ is the vector of initially unknown data belonging to the uncertainty set $\Z \subseteq \R^{N_z}$. Solving problem \eqref{eq:sec1:two.stage} is difficult in general, since $\Z$ might include an infinite number of scenarios, and hence different values of $\y$ might be optimal for different realizations of $\z$. In fact, finding optimal $\x$ is an NP-hard problem \citep{guslitser2002uncertainty}. To address this difficulty, several approaches have been proposed. The first one is to use so-called decision rules which explicitly formulate the second-stage decision $\y$ as a function of $\z$, and hence the function parameters become first-stage decisions next to $\x$; see \cite{ben2004adjustable}. Another approach is to partition $\Z$ into subsets and to assign a separate copy of $\y$ to each of the subsets. The partitioning is then iteratively refined, and the decisions become increasingly \emph{customized} to the outcomes of $\z$.

In this paper, we consider a third approach to \eqref{eq:sec1:two.stage} known as $K$-adaptability. There, at most $K$ possible wait-and-see decisions $\y_1, \ldots, \y_K$ are allowed to be constructed, and the decision maker must select one of those. The values of the possible $\y_k$'s become the first-stage variables, and the problem boils down to
\begin{equation}
    \min_{\x \in \X, \y \in \Y^K} \max_{\z \in \Z} \; \min_{k \in \K} \big \{\bc(\z)^\intercal \x + \d(\z)^\intercal \y_k: \T(\z)\x + \W(\z)\y_k \leq \h(\z)\big \},
    \label{eq:sec1:k_adaptability_problem}
\end{equation}
where $\K = \{1, \ldots, K\}$ and $\Y^K = \bigtimes_{k = 1}^K {\Y}$. Although the solution space of \eqref{eq:sec1:k_adaptability_problem} is finite-dimensional, it remains an NP-hard problem. For certain cases, \eqref{eq:sec1:k_adaptability_problem} can be equivalently rewritten as a mixed integer linear programming (MILP) model \citep{hanasusanto2015k}. 

The above formulation requires that for given $\x \in \X$ and $\z \in \Z$, there is at least one decision $\y_k$, $k \in \K$ satisfying $\T(\z)\x + \W(\z)\y_k \leq \h(\z)$, and among those one (or more) minimizing the objective. Looking at \eqref{eq:sec1:k_adaptability_problem} from the point of view $\y_k$, we can say that for each $\y_k$, we can identify a subset $\Z_k$ of $\Z$ for which a given $\y_k$ is optimal {among $K$ selected recourses}. The union of sets $\Z_k$, $k \in \K$ is equal to $\Z$ although they need not be mutually disjoint (but a mutually disjoint partition of $\Z$ can be constructed). Consequently, solving \eqref{eq:sec1:k_adaptability_problem} involves implicitly (i) clustering $\Z$, and (ii) optimizing the per-cluster decision so that the objective function corresponding to \emph{the most difficult cluster} is minimized. Such simultaneous clustering and per-cluster optimization also occur, for example in retail. A line of $K$ products is to be designed to attract the largest possible group of customers. The customers are clustered into $K$ groups, and the nature of the products is guided by the cluster characteristics.

In this manuscript, we focus on the general {MILP} $K$-adaptability case for which the only existing solution approach is the $K$-adaptability branch-and-bound ($K$-B\&B) algorithm of \cite{subramanyam2020k}. {Other methods to solve the $K$-adaptability problem have been proposed by \cite{hanasusanto2015k} that deals with binary decisions, and \cite{ghahtarani2023double} that assumes integer first-stage decisions}. The $K$-B\&B algorithm, as opposed to the top-down partitioning of $\Z$ of \cite{bertsimas2016multistage} or \cite{postek2016multistage}, proceeds by gradually building up discrete subsets $\barZ_k$ of scenarios. In most practical cases, a solution to \eqref{eq:sec1:k_adaptability_problem}, where $\y_1, \ldots, \y_K$ are feasible for large $\barZ_1, \ldots, \barZ_K$, is also feasible to the original problem. The problem, however, lies in knowing which scenarios should be grouped together. In other words, a decision needs to be made on which scenarios of $\Z$ should be responded to with the same decision. How well this question is answered, determines the (sub)optimality of $\y_1, \ldots, \y_K$. In \cite{subramanyam2020k}, a search tree is used to determine the best collection (see Section \ref{sec:background} for details). However, this approach suffers from exponential growth.

We introduce a method for learning the best strategy to explore this tree. In particular, we learn which nodes to evaluate next in depth-first search \emph{dives} to obtain good solutions faster. These predictions are made using a supervised machine learning (ML) model. {As the input does not range for different instance sizes, or values of $K$, as will be explained in future sections, the ML model does not require difficult architectures. Standard ML models, such as feed-forward neural networks, random forests, or support vector machines can be used for this work.}

Due to the supervised nature, some \emph{oracle} is required to be imitated. In the design of this oracle, we are partly inspired by Monte Carlo tree search (MCTS) \citep{browne2012survey}, which is often used for exploring large trees. Namely, the training data is obtained by exploring $K$-B\&B trees via an adaptation of MCTS (see Section \ref{sec:tree_engineering}). The scores given to the nodes in the MCTS-like exploration are stored and used as labels in our training data.

In the field of solving MILPs, learning node selection to speed up exploring the B\&B tree has been done, \emph{e.g.}, by \cite{he2014learning}. Here, a node selection policy is designed by imitating an oracle. This oracle is constructed using the optimal solutions of various MILP data sets. More recently, \cite{khalil2022mip} used a graph neural network to learn node selection. {Our method distinguishes itself from these approaches as we specifically use the nature of our problem. Namely, in the design of the node selection strategy, we use the actual meaning of selecting a node; adding a scenario to a subset. Therefore, we try to learn what characteristics (or features) scenarios that should belong to the same subset have.}

For an overview on ML for learning branching policies in B\&B, see \cite{bengio2020machine}. There has also been done a vast amount of research on applying MCTS directly to solving combinatorial problems. In \cite{sabharwal2012guiding} a special case of MCTS called Upper Confidence bounds for Trees (UCT), is used for designing a node selection strategy to explore B\&B trees (for MIPs). In \cite{khalil2022finding} MCTS is used to find the best backdoor (\textit{i.e.}, a subset of variables used for branching) for solving MIPs. \cite{loth2013bandit} have used MCTS for enhancing constraint programming solvers, which naturally use a search tree for solving combinatorial problems. For an elaborate overview on modifications and applications of MCTS, we refer to \cite{swiechowski2022monte}.

The remainder of the paper is structured as follows. In Section \ref{sec:background} we describe the inner workings of the $K$-adaptability branch-and-bound to set the stage for our contribution. In Section \ref{sec:methodology} we outline our ML methodology along with the data generation procedure. Section \ref{sec:experiments} discusses the results of a numerical study, and Section \ref{sec:conclusion} concludes with some remarks on future works. 

\section{Preliminaries}  \label{sec:background}

\begin{figure}[htbp!]
    \centering
    \includegraphics{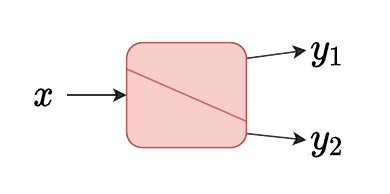}
    \caption{A framework of the $K$-adaptability problem, where we split the uncertainty set (red box) in $K=2$ parts. Here, $\x$ represents the first-stage decisions, and $\y_1$ with $\y_2$ those of the second-stage.}
    \label{fig:k_adapt_example}
\end{figure}

It is instructive to conceptualize a solution to \eqref{eq:sec1:k_adaptability_problem} as a solution to a nested clustering and optimization-for-clusters methodology. As already mentioned in Section \ref{sec:introduction}, a feasible solution to \eqref{eq:sec1:k_adaptability_problem} can be used to construct a partition of the uncertainty set into subsets $\Z_1, \ldots, \Z_K$ such that $\bigcup_{k=1}^K \Z_k = \Z$. Here, decision $\y_k$ is applied in the second stage if $\z \in \Z_k$. The decision framework associated with a given solution is illustrated in Figure \ref{fig:k_adapt_example}. 

\noindent For such a \emph{fixed partition} the corresponding optimization problem becomes
\begin{align}
    \min_{\x \in \X, \y \in \Y^K} &  \max_{k \in \K} \max_{\z \in \Z_k} \left\{ \bc(\z)^\intercal \x + \d(\z)^\intercal \y_k \right\} \label{sec2:per_partition_problem} \\
    \text{s.t.} \quad & \T(\z)\x + \W(\z)\y_k \leq \h(\z). \nonumber 
\end{align}
The optimal solution to \eqref{eq:sec1:k_adaptability_problem} also corresponds to an optimal partitioning of $\Z$, and the optimal decisions of \eqref{sec2:per_partition_problem} with that partitioning. Finding an optimal partition and the corresponding decisions has been shown to be NP-hard by \cite{bertsimas2010finite}. For that reason, \cite{subramanyam2020k} have proposed the $K$-B\&B algorithm. There, the idea is to gradually build up a collection of finite subsets $\barZ_1, \ldots, \barZ_K$, such that for each $k \in \K$ an optimal solution to \eqref{sec2:per_partition_problem} with $\Z_k = \barZ_k$ is also an optimal solution to \eqref{eq:sec1:k_adaptability_problem}.

\begin{figure}[htbp!]
    \centering
    \includegraphics[width=0.65\columnwidth]{
    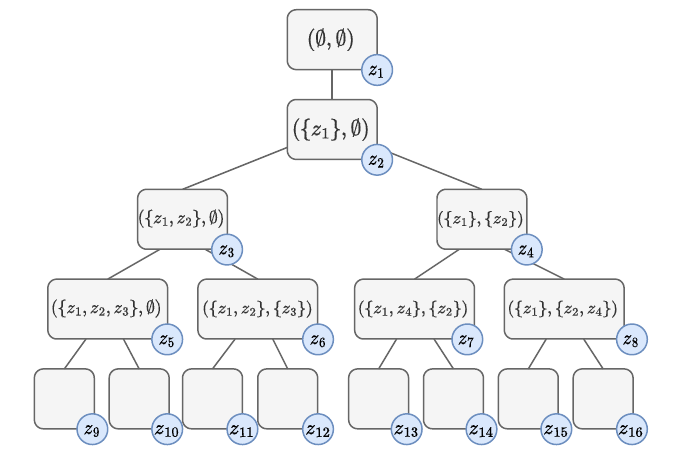}
    \caption{Search tree for $K$-adaptability branch-and-bound ($K = 2$).}
    \label{fig:tree}
\end{figure}

The algorithm follows a master-subproblem approach. The master problem solves \eqref{eq:sec1:k_adaptability_problem} with $K$ finite subsets of scenarios. The subproblem finds the scenario for which the current master solution is not robust. The number $K$ of possible assignments of this new scenario to one of the existing subsets gives rise to using a search tree. Each tree node corresponds to a partition of all scenarios found so far into $K$ subsets. The goal is to find the node with the best partition. An illustration of the search tree is given in Figure \ref{fig:tree}. The tree grows exponentially and thus only (very) small-scale problems can be solved in reasonable time. The method we propose in the next section learns a good node selection strategy with the goal of converging to the optimal solution much faster than $K$-B\&B. 

\paragraph{\textit{Master problem}.}
This problem solves the $K$-adaptability problem \eqref{eq:sec1:k_adaptability_problem} with respect to the currently found scenarios grouped into $\barZ_k \subset \Z$ for all $k \in \K$. For a collection $\barZ_1, \ldots, \barZ_K$, the problem formulation is defined as follows:
\begin{align}
\min_{\theta \in \R, \x\in \X, \y \in \Y^K} & \theta \label{eq:master_problem} \\
\text{s.t.} \qquad & \bc(\z)^\intercal \x + \d(\z)^\intercal \y_k \leq \theta, & \forall \z \in \barZ_k, \forall k \in \K, \nonumber \\
 & \T(\z)\x + \W(\z)\y_k \leq \h(\z), \quad \quad \quad & \forall \z \in \barZ_k, \forall k \in \K, \nonumber
\end{align}
where $\theta$ is the current estimate of the objective function value. We denote the optimal solution of \eqref{eq:master_problem} by the triplet $(\theta^*, \x^*, \y^*)$. 

\paragraph{\textit{Subproblem}.} The subproblem aims to find a scenario $\z$ for which the current master solution is infeasible. That is, a scenario is found such that for each $k$, at least one of the following is true:
\begin{itemize}
    \item the current estimate of $\theta^*$ is too low, \textit{i.e.}, $\bc(\z)^\intercal \x^\ast + \d(\z)^\intercal \y^\ast_k > \theta^\ast$;
    \item at least one of the original constraints is violated, \textit{i.e.}, $\T(\z)\x^\ast + \W(\z)\y_k^\ast > \h(\z)$.
\end{itemize}
If no such scenario exists, we define the solution $(\theta^*, \x^*, \y^*)$ as a robust solution. When such a scenario $\z^*$ does exist, the solution is not robust and the newly-found scenario is assigned to one of the sets $\barZ_1, \ldots, \barZ_K$.

\begin{definition}
A solution $(\theta^*, \x^*, \y_1^*, \ldots, \y_K^*)$ to \eqref{eq:master_problem} is robust if
$$
\forall \z \in \Z, \exists k\in \K : \ \T(\z)\x^\ast + \W(\z)\y_k^\ast \leq \h(\z), \ \bc(\z)^\intercal \x^\ast + \d(\z)^\intercal \y_k^\ast \leq \theta^\ast.
$$
\end{definition}

\begin{example}
Consider the following master problem
\begin{align*}
    \min_{\theta, \x, \y} \quad & \theta\\
    \text{s.t.} \quad & \theta \in \R, \x\in \{0, 1\}^2, \y \in \{0, 1\}^2,\\
    & \z^\intercal \x + [2z_1, 2z_2] \y_k \leq \theta, & \forall \z \in \barZ_k, \forall k \in \K, \\
    & \y_k \geq \1 - \z, & \forall \z \in \barZ_k, \forall k \in \K, \\
    & \x + \y_k \geq \1, & \forall k \in \K,
\end{align*}
where $\z \in \{-1, 0, 1\}^2$ and $\K = \{1, 2\}$. We look at three solutions for different groupings of $\z_1 = [1, 1]^\intercal$, $\z_2 = [0, 1]^\intercal$, and $\z_3 = [0, 0]^\intercal$. One of them is not robust, and the other ones are but have a different objective value due to the partition. The first partition we consider is $\barZ_1 = \{\z_1\}, \ \barZ_2 = \{\z_2\}$. The corresponding solution is $\theta^* = 2$, $\x^* = [1, 1]^\intercal$, $\y^*_1 = [0, 0]^\intercal$, and $\y^*_2 = [1, 0]^\intercal$. This solution is not robust, since for $\z_3 = [0, 0]^\intercal$ the following constraint is violated for all $k \in \K$:
\begin{equation*}
        \y_1 \ngeq \1 - \z_3 \iff     \begin{bmatrix} 
                                           0 \\
                                           0 \\
                                        \end{bmatrix} 
                                        \ngeq 
                                        \begin{bmatrix} 
                                           1 \\
                                           1 \\
                                        \end{bmatrix} 
                                        - 
                                        \begin{bmatrix} 
                                           0 \\
                                           0 \\
                                        \end{bmatrix},
        \qquad
        \y_2 \ngeq \1 - \z_3  \iff     \begin{bmatrix} 
                                           1 \\
                                           0 \\
                                        \end{bmatrix} 
                                        \ngeq 
                                        \begin{bmatrix} 
                                           1 \\
                                           1 \\
                                        \end{bmatrix} 
                                        - 
                                        \begin{bmatrix} 
                                           0 \\
                                           0 \\
                                        \end{bmatrix}. \\
\end{equation*}
Next, consider $\barZ_1 = \{\z_1, \z_2\}, \ \barZ_2 = \{\z_3\}$, with solution $\theta^* = 3$, $\x^* = [0, 1]^\intercal$, $\y^*_1 = [1, 0]^\intercal$, $\y^*_2 = [1, 1]^\intercal$. There is no $\z$ that violates this solution, which makes it robust. The third partition is $\barZ_1 = \{\z_1, \z_3\}, \ \barZ_2 = \{\z_2\}$, with solution $\theta^* = 4$, $\x^* = [0, 0]^\intercal$, $\y^*_1 = [1, 1]^\intercal$, and $\y^*_2 = [1, 1]^\intercal$. This solution is also robust. Their objective values are 3 and 4, which shows that different partitions can significantly influence solution quality.
\end{example}

Mathematically, we formulate the subproblem using the big-$M$ reformulation:
\begin{align}
\max_{\zeta, \z, \gamma} \quad & \zeta  \label{eq:sub_problem}\\
\text{s.t.} \quad & \zeta \in \R, \quad \z \in \Z, \quad \gamma_{kl} \in \{0, 1\}, \quad (k, l) \in \K \times \L, \nonumber \\
& \sum_{l \in \L} \gamma_{kl} = 1, & \forall k \in \K,  \nonumber \\
    & \zeta + M(\gamma_{k0}-1) \leq \bc(\z)^\intercal \x^* + \d(\z)^\intercal \y^*_k - \theta^*, \qquad & \forall k \in \K, \nonumber \\
    & \zeta + M(\gamma_{kl}-1) \leq \boldsymbol{t}_l(\z)^\intercal \x^* + \boldsymbol{w}_l(\z)^\intercal \y^*_k - h_l(\z), & \forall l \in \L, \forall k \in \K, \nonumber
\end{align}
where $M$ is some big scalar and $\L$ is the index set of constraints. When $\zeta \leq 0$, we have not found any violating scenario, and the master solution is robust. Otherwise, the found $\z^*$ is added to one of $\barZ_1, \ldots, \barZ_K$. Then, the master problem is re-solved, so that a new solution ($\theta^*$, $\x^*$, $\y^*$), which is guaranteed to be feasible for $\z^*$ as well, is found.

Throughout the process, the key issue is to which of $\barZ_1, \ldots, \barZ_K$ to assign $\z^*$. As this cannot be determined in advance, all options have to be considered. Figure \ref{fig:tree} illustrates that from each tree node we create $K$ child nodes, each corresponding to adding $\z^\ast$ to the $k$-th subset: 
\begin{equation*}
   \tau^k = \{\barZ_1, \ldots, \barZ_k \cup \{\z^*\}, \ldots, \barZ_K\}, \quad \forall k \in \K,
\end{equation*}
where $\tau^k$ is the partition corresponding to the $k$-th child node of node $\tau$. If a node is found such that $\theta^*$ is greater than or equal to the value of another robust solution, then this branch can be pruned as the objective value of the master problem cannot improve if more scenarios are added.

The pseudocode of $K$-B\&B is given in Algorithm~\ref{algorithm:k_bb} in the Appendix. The tree becomes very deep for 2SRO problems where many scenarios are needed for a robust solution. Moreover, the tree becomes wider when $K$ increases. Due to these issues, solving the problem to optimality becomes computationally intractable in general. Our goal is to investigate if ML can be used to make informed decisions regarding the assignment of newly found $\z^*$ to the discrete subsets, so that a smaller search tree has to be explored in a shorter time before a high-quality robust solution is found.

\section{ML methodology} \label{sec:methodology}
We propose a method that enhances the node selection strategy of the $K$-B\&B algorithm. It relies on four steps:
\begin{enumerate}
    \item \textit{\textbf{Decision on what and how we want to predict:}} Section \ref{sec:node_selection}.
    \item \textit{\textbf{Feature engineering:}} Section \ref{sec:features}.
    \item \textit{\textbf{Label construction:}} Section \ref{sec:learn_strategy}.
    \item \textit{\textbf{Using partial $K$-B\&B trees for training data generation:}} Section \ref{sec:tree_engineering}.
\end{enumerate}
All these steps are combined into the \textsc{$K$-B\&B-NodeSelection} algorithm (Section \ref{sec:algorithm}).

\subsection{Learning setup} \label{sec:node_selection}
\begin{figure}[htbp!]
    \centering
        \includegraphics[width=0.65\columnwidth]{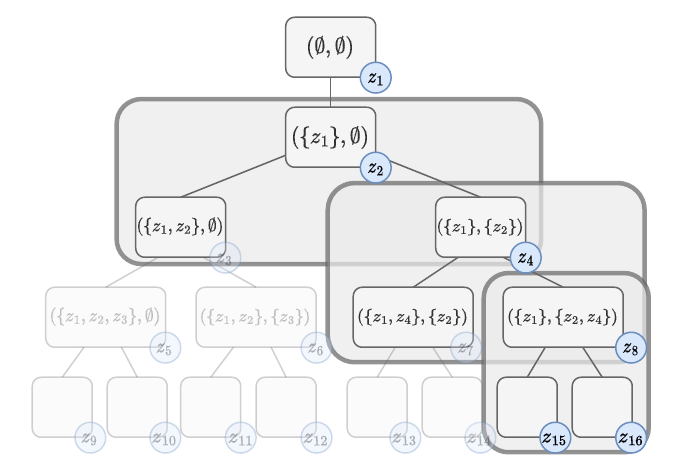}
    \caption{The scope of node selection}
    \label{fig:node_select}
\end{figure}

As there is no clearly well-performing node selection strategy for $K$-B\&B, we cannot simply try to imitate one. Instead, we investigate what choices a good strategy would make. We will focus on learning how to make informed decisions about the order of inspecting the children of a given node. The scope of our approach is illustrated with rounded-square boxes in Figure \ref{fig:node_select}.

To rank the child nodes in the order they ideally be explored, one of the ways is to have a certain form of child node information whether selecting this node is \emph{good} or \emph{bad}. In other words, how likely is a node to guide us towards a high-quality robust solution fast. Indeed, this shall be exactly the quantity that we predict with our model. To train such a model, we will construct a dataset consisting of the following input-output pairs:
\begin{itemize}
    \item \textbf{Input.} feature vector $F$ of the decision to insert a scenario to a subset (Section \ref{sec:features})
    \item \textbf{Output.} $[0, 1]$ label that informs how good a given insertion was, based on an ex-post constructed strategy; `what would have been the best node selection strategy, had we known the entire tree?' (Section \ref{sec:learn_strategy}).
\end{itemize}

\begin{figure}[htbp!]
    \centering
        \includegraphics[width=0.9\columnwidth]{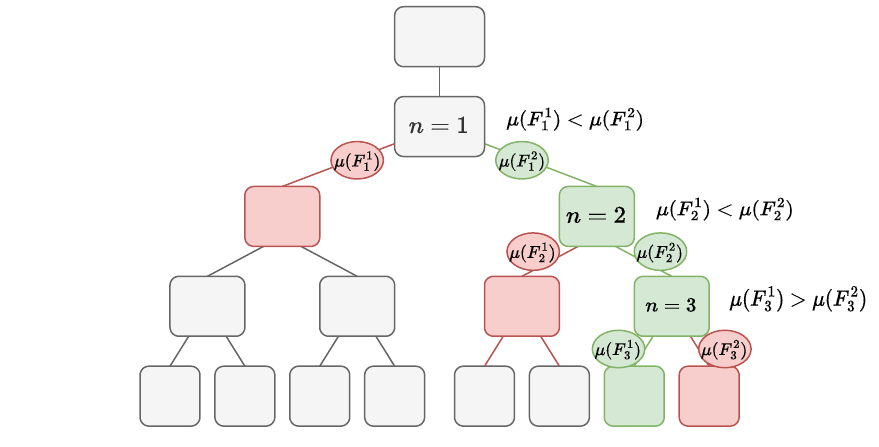}
    \caption{Node selection with ML predictions. Node selections are decided by the prediction of the function $\mu$ with input features $F_n^k$, where $n$ is the node from which a selection is made and $k$ relates to its $k$-th child node.}
    \label{fig:node_select_ml}
\end{figure}

We gather this data by creating a proxy for an oracle (see Section \ref{sec:tree_engineering}). Once the predictive model is trained, we can apply it to the search tree where in each iteration we predict for all $K$ child nodes a score $\mu$, in the interval $[0,1]$. The child node with the highest score is explored first (see Figure \ref{fig:node_select_ml}).
\noindent Formally, the working order of the trained method shall be as follows:
\begin{enumerate}
    \item \textbf{Process node.     } Solve master and subproblem in current node $n$. If prune conditions hold, then select a new node. Otherwise, continue to Step 2.
    \item \textbf{Compute features.   } For each one of the $K$ child nodes, generate a vector of features $F^1_n, \ldots, F^K_n$ (see Section \ref{sec:features} for details).\label{alg_node_select:features}
    \item \textbf{Predict.  } For each child node $k$, predict the goodness score $\mu (F^k_n)$.
    \item \textbf{Node selection.   } Select the child node with the highest score.
\end{enumerate}
It is desirable for an ML tool to be applicable to data that is different from the one which it is trained on. In the context of an ML model constructed to solve optimization problems, this means the potential to use a given trained method on various optimization problems of various sizes, with different values of $K$. Indeed, we shall demonstrate the generality of our method. First, we note that each parent-child node pair corresponds to one data point, and hence, training on a problem with a certain $K$ does not prevent us from using it for different $K$. Next, we construct our training dataset so that the model becomes independent of (i) the instance size, and (ii) the type of objective function and constraints. This will be explained in Section \ref{sec:features}.

\subsection{Feature engineering} \label{sec:features}
If we take a look at a single rounded box in Figure \ref{fig:node_select}, the ML model we are about to train is going to give a goodness score which will depend on the parent node and the scenario. Therefore, we need to design features which we group into (i) state features that describe the master problem and the subproblem solved in the parent node, (ii) scenario features that describe the assignment of a newly-found scenarios to one of the subsets $\barZ_k$. In what follows, we present the feature list:
\begin{enumerate}
    \item \textbf{State features.} This input describes the parent node $n$, \textit{i.e.}, the current state of the algorithm. Different states might benefit from different strategies. Hence, information on the current node might increase the prediction performance. This also means that all the child nodes have the same state features: $\s_n^k = \s_n$ for all $k \in \K$. To scale the features, we always initialize a tree search with a so-called \emph{initial run}. This is a dive with random child-node selections, where we stop until robustness is reached. The following values are taken from this initial run: 
    \begin{itemize}
        \item $\theta^0$: the objective function value of the robust solution found in the initial run.
        \item $\zeta^0$: the violation of the root node. 
        \item $\kappa^0$: depth reached in the initial run.
    \end{itemize}
    In the experiments, multiple initial runs are done. The averages of $\theta^0$, $\zeta^0$, and $\kappa^0$ over the dives are then used for scaling. {Additional meta-features for $K$ and the dimensions of $\Z$, $\X$, and $\Y$ could be added. We omitted these as in our experimental setup we do not mix training data for different combinations of these values.}
    \item \textbf{Scenario features.} Intuitively, scenarios contained in the same group should have similar characteristics. Therefore, the features for each node are constructed in the following way: each newly found scenario $\z^*$ is assigned a set of characteristics, or attributes. Based on the attributes of the new scenario, and the attributes of the scenarios already grouped into the $K$ subsets, we formulate the input of one data point. Some of the scenario attributes can be directly determined from the master problem. Others are extracted from easily solvable optimization problems: the \textit{deterministic problem} and the \textit{static problem}. The deterministic problem is a version of the problem where $\z^*$ is the only scenario. {The solution to this problem gives information on the optimal objective and decisions found, in the most naive sense. One would expect that for good-performing subset formations, the optimal decisions of its scenarios would have similarities}. For the static problem, we first solve for a single $\y$ (no adaptability) for all $\z \in \Z$. Then, for the obtained $\x^*$ and for a given $\z^*$, we solve for the best $\y$. {The solution to this problem gives additional information as it has a sense of embedded robustness.}
\end{enumerate}

As our goal is to use the model also on different problems, we engineer the features to keep them independent from the problem size or type. In Tables \ref{tab:state_features}-\ref{tab:attributes}, we give an outline of the two feature types that we construct. For a detailed description of how they are computed, we refer the reader to Appendix~\ref{app:attributes}.

\begin{table}[htbp!]
	\caption{State features. $\theta^0$, $\zeta^0$, and $\kappa^0$ are as defined above. $\theta^p$ is the objective value of the parent node, $\zeta^p$ is the violation of the parent node, and $\kappa$ is the depth of the current node.}
	\footnotesize
	\centering
	\begin{tabular}{llp{60mm}c}
        \toprule
		 \textbf{Num.} & \textbf{State feature name} & \textbf{Description} & \textbf{Calculation}\\
		\midrule
		1 & Objective & Relative objective value of this node to the first robust solution & $\theta/\theta^0$    \\ \\
		2 & Objective difference & Ratio of objective to that of the parent node & $\theta/\theta^p$ \\ \\
        3 & Violation & Relative violation with respect to the first violation found & $\zeta/\zeta^0$ \\ \\
        4 & Violation difference & Ratio of violation to that of the parent node & $\zeta/\zeta^p$ \\ \\
        5 & Depth & Relative depth of this node to the depth of the fist robust solution & $\kappa/\kappa^0$ \\
		\hline
	\end{tabular}
	\label{tab:state_features}
\end{table}

\begin{table}[htbp!]
	\caption{Attributes assigned to scenario $\z^*$.} 
	\footnotesize
	\centering
	\begin{tabular}{l|llp{60mm}}
		\toprule
		 {} & \textbf{Num.} & \textbf{Attribute name} & \textbf{Description}\\
		\midrule
		 {} & 1 & Scenario values & Vector of scenario values $\z^*$  \\
		\midrule
		\multirowcell{2}{\textit{Master} \\ \textit{problem}} & 2 & Constraint distance     & A measure for the change of the feasibility region when $z^*$ is added to a subset. We look at the distance between the constraints already in the master problem, and the one to be added.    \\
		& 3 & Scenario distance     &  With this attribute we measure how far away $\z^*$ is from \textit{not} being a violating scenario, for each of the $k$ subsets. This is done by looking at the constraints in the space of $\Z$, given the current solutions $\x$ and $\y^k$.\\ 
		& 4 & Constraint slacks   &  The slack values of the uncertain constraint per subset decisions. \\
		\midrule
		\multirowcell{2}{\textit{Deterministic} \\ \textit{problem}} &  5 & Objective value & The objective value of the deterministic problem \\
		& 6 & First-stage decisions & First-stage decisions of deterministic problem \\
		& 7 & Second-stage decisions & Second-stage decisions of deterministic problem \\
		\midrule
		 \multirowcell{2}{\textit{Static} \\ \textit{problem}} & 8 &  Objective value &  The objective value of the static problem  \\
		& 9 &  second-stage decisions & Second-stage decisions of static problem \\
		\bottomrule
	\end{tabular}
	\label{tab:attributes}
\end{table}
The state features in Table \ref{tab:state_features} are readily problem-independent. However, this is not the case for the scenario attributes in Table \ref{tab:attributes}. They are, for example, dependent on the instance size. {To deal with the size dependency, and to adopt the actual meaning of node selection (\emph{i.e.}, placing a scenario to a group of other scenarios), we introduce the \emph{attribute distance} as a feature to our model. This is the distance between scenario-subset pairs in terms of the attributes attached to the scenarios, as a proxy to the unknown implications the addition of the scenario would have on the solution of the problem. For each attribute, the scenario-subset distance is taken by the Euclidean distance from the attribute of the new scenario to the average of the attribute of the scenarios already in the subset.} This feature is described as: 
\begin{equation}
    \delta^k_f = \frac{\lVert\a^{k, \z^*}_f - \frac{1}{|\barZ_k|}\sum_{\z \in \barZ_k}{\a_f^{k, \z}}\rVert_2}{{\text{length}}(\a^{k, \z^*}_f)},  \qquad \forall k \in \K, \forall f \in \{1, \ldots, 9\},  \label{eq:att_dist}
\end{equation}
where $\delta_f^k$ is the attribute distance of the new scenario $\z^*$ to subset $\barZ_k$ and $\a_f^{k, \z}$ is the data vector related to the $f$-th attribute (of Table \ref{tab:attributes}) for the $k$-th child node. {The attribute distance is scaled by the length of the attribute vector, denoted by `length' in the denominator of \eqref{eq:att_dist}, to control for varying attribute vectors in the feature value. Larger attribute vectors would otherwise result in higher feature values.} Then, the \emph{scenario feature} vector of the $k$-th child node is defined as $\d^k_n = [\delta^k_1, \ldots, \delta^k_9]^\intercal$. 

Then, the input of the ML model for the $k$-th child of node $n$ is given by $F^k_n = [\s_n \quad \d_n^k]^\intercal$. 

\begin{figure}[htbp!]
    \centering
        \includegraphics[width=\columnwidth]{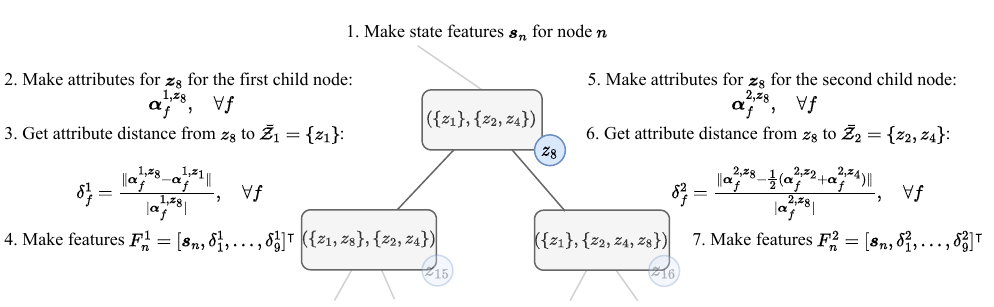}
    \caption{Example of a feature generation procedure for node $n$ and its two child nodes.}
    \label{fig:features_example}
\end{figure}

In Figure \ref{fig:features_example} we outline the feature generation procedure. Steps 2 and 5 indicate how new attributes need to be generated for every child node. In practice many attributes are the same for all subsets, and the attributes that do differ are the master problem-based ones, which are easily computed.

\subsection{Label construction} \label{sec:learn_strategy}
To learn how to assign a new scenario to one of the subsets $\barZ_k$, we need another piece of the input-output pairs in our database -- labels that would indicate how good, ex post, it was to perform a given assignment, \textit{i.e.}, how likely a given assignment is to lead the search strategy towards a \emph{good solution}. We shall assume that given a $K$-B\&B tree, a good solution is a robust solution with objective that belongs to the best $\alpha \%$ of the found robust solutions. We will now introduce a notion of scenario-to-set assignments and illustrate a method of constructing labels $q$.

We define $p_\nu$ -- the probability of node $\nu$ leading to a good robust solution.  If this node is selected, one of a finite, possibly large, number of \emph{leaves} (\textit{i.e.}, terminal nodes) is reached with depth-first search. Then, taking $g_\nu$ to be the number of good solutions and $M_\nu$ as the number of leaves under node $\nu$, we define $p_\nu = g_\nu / M_\nu$ as the fraction of leaves with a good robust solution. We define a node $\nu$ as good if the probability $p_\nu$ of it leading to a good robust solution is higher than some threshold $\epsilon$. This is computed by the quality value $q_\nu = \1_{\{p_\nu \geq \epsilon\}}$.

\begin{figure}[htbp!]
    \centering
        \includegraphics[width=0.6\columnwidth]{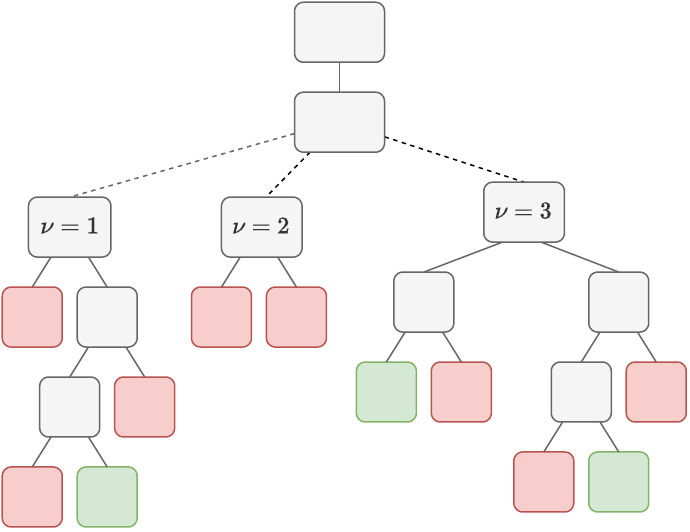}
    \caption{Example of a tree where three nodes are considered. The coloured nodes are leaves; green and red nodes are good and bad robust solutions, respectively.}
    \label{fig:tree_suc_pred}
\end{figure}

\begin{example}
Consider the tree in Figure \ref{fig:tree_suc_pred}. For three nodes in the search tree, its subtree (consisting of all its successors) and leaves (coloured nodes) are shown. Red coloured nodes represent bad and green nodes represent good robust solutions. Then, if we pick $\epsilon = \frac{1}{5}$ as threshold, the corresponding  success probabilities $p_\nu$ and the quality values of the three nodes are:
\begin{equation*}
    p_1 = \frac{1}{4}, \
    p_2 =  0, \ 
    p_3 = \frac{2}{5}, \
    q_1 = \1_{\{p_1 \geq \frac{1}{5}\}} = 1, \ 
    q_2 = \1_{\{p_2 \geq \frac{1}{5}\}} = 0, \
    q_3 = \1_{\{p_3 \geq \frac{1}{5}\}} = 1.
\end{equation*}
The first and third nodes would have been good node selections, whereas the second selection would have been a bad one.
\end{example}

In practice, we do not know for a node $\nu$ how many underlying leaves are good and bad. This is the reason that in this method we will predict $q_\nu$ with a model $\mu$. We will also call this model the strategy model (or function) since it guides us in making node selections.

\subsection{Training data generation} \label{sec:tree_engineering}
Our goal is to learn a supervised ML model $\mu(F_\nu) = \hq_\nu$,
where $\mu$ is the strategy function, $F_\nu$ the input features, and $\hq_\nu$ the prediction of the quality of moving to node $\nu$. To train this model, we first need to generate training data. Given a single tree, the difficulty of generating data does not lie in the input, but in the output: an expert is needed to determine the correct values of $\hq_\nu$ for its nodes. Consider the following; if nodes of the entire tree are \emph{processed} (\textit{i.e.}, solving the master problem and the subproblem), we could easily find the paths from the root to good solutions. Then, all the nodes in these paths would get the values $\hq_\nu = 1$, and all others $\hq_\nu = 0$. Or equivalently, the success probabilities would be set to $p_\nu > \epsilon$. Since the trees grow exponentially, this approach is not practical.

\begin{figure}[htbp!]
    \centering
        \includegraphics[width=\columnwidth]{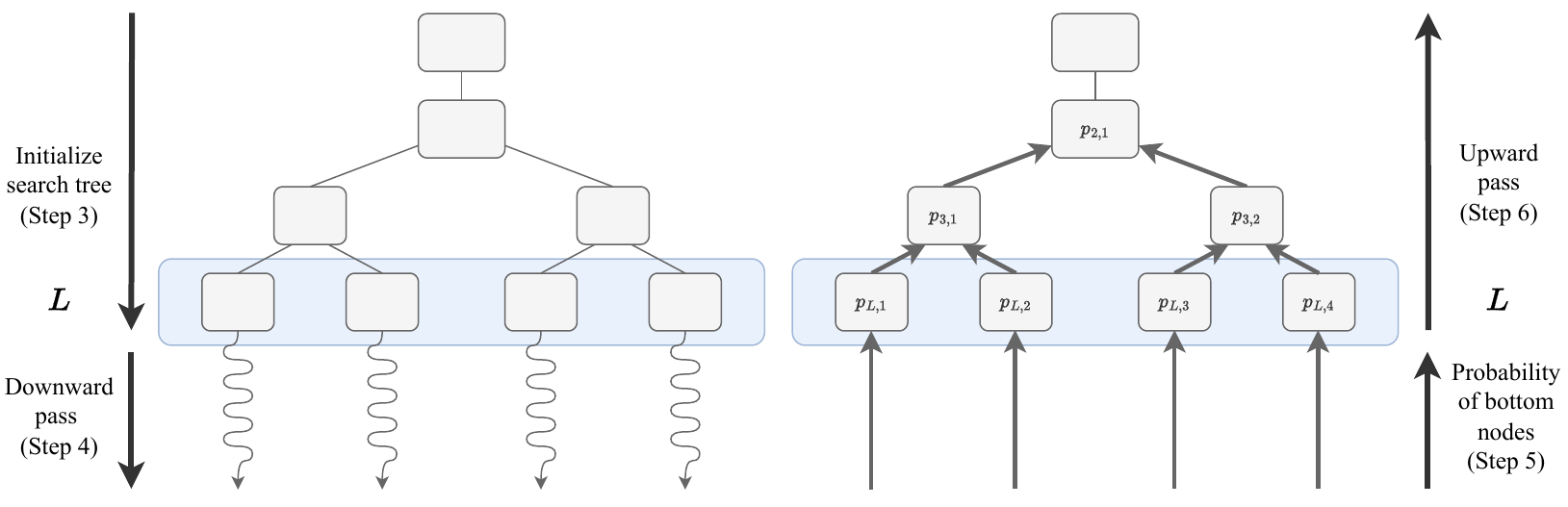}
    \caption{Downward and upward pass for generating training data. The blue layer represents the $L$-th layer in the tree from which random runs are made.}
    \label{fig:training}
\end{figure}

As already mentioned in the introduction, a popular method for exploring intractable decision trees is Monte Carlo tree search (MCTS) \citep{browne2012survey}. By randomly running deep in the tree, we can gather information on the search space. In our method, we use the idea of random runs to mimic an expert, and hence, label the data points. Generating training data is done as follows (see Figure \ref{fig:training}):
\begin{enumerate}
    \item \textbf{Get instance. } Generate an instance of a 2SRO problem. \label{train_step:get_instance}
    \item{\textbf{Initial run.   }} One (or multiple) dives are executed to gather feature information.
    \item \textbf{Initialize search tree.  } Process all nodes up to a predetermined level $L$ of the tree. Generate features for each of the explored nodes.
    \item \label{train_step:downward}\textbf{Downward pass.    } Per node $\nu \in \{1, \ldots, N_L\}$ of the $L$-th layer ($N_l$ is the number of nodes of the $l$-th layer), perform dives for a total of $R$ times.
    \item \textbf{Probability of bottom nodes.  } Set probability $p_{L,\nu}$ of each node $\nu \in \{1, \ldots, N_L\}$ in layer $L$ as $p_{L,\nu} = \frac{g_{\nu}}{R}$ where $g_{\nu} \in \{0, \ldots,  R\}$ is the number of good solutions from the samples of the $\nu$-th node in layer $L$.
    \item \textbf{Upward pass.    } Propagate the probabilities $p_{l,\nu}$ upwards through the tree, for all nodes $\nu \in \{1, \ldots, N_l\}$, for all levels $l \in \{L-1, \ldots, 2\}$, as follows:
    \begin{align*}
        p_{l,\nu} & = \mathbb{P}(\text{at least one child node is successful}) \\
        & = 1 - \mathbb{P}(\text{no successful child nodes}) \\
        & = 1 - \prod_{k \in \K} (1-p^k_{l,\nu}),
    \end{align*} 
    where $p^k_{l,\nu}$ is the $k$-th child of node $\nu$ of layer $l$ (while this child node is in the ($l+1$)-th level). Note that $p^k_{l,\nu} = p_{l+1,\nu'}$ for some $\nu \in \{1, \ldots, N_l\}$ and $\nu' \in \{1, \ldots, N_{l+1}\}$.
    \item \textbf{Label nodes.  } Determine the label $\hq_{l, \nu}$ for nodes $\nu \in \{1,\ldots, N_l\}$ for levels $l \in \{2, \ldots, L\}$. \label{train_step:upward_pass}
\end{enumerate}

The advantage of this structure is that both bad and good decisions are well represented in the dataset. However, it only consists of input-output pairs of the top $L$ levels. This is not necessarily a disadvantage as good decisions at start can be expected to be more important. The above generation method is applied per instance, thus it can be parallelized, like Step \ref{train_step:downward}. 

\subsection{Complete node selection algorithm} \label{sec:algorithm}
We now combine all the steps above into one algorithm which is, essentially, a variant of $K$-B\&B enhanced with: (i) node selection, (ii) feature engineering, and (iii) training data generation (see Algorithm \ref{algorithm:k_bb_node_selection}). Our \textsc{$K$-B\&B-NodeSelection} algorithm has two preprocessing steps: 
\begin{enumerate}
    \item \textsc{StrategyModel} (Procedure~\ref{proc:strategy_model} in Appendix~\ref{app:pseudocodes}): Generate the data applying \ref{train_step:get_instance}-\ref{train_step:downward}, and train the ML model.
    \item \textsc{InitialRun} (Procedure~\ref{proc:init_run} in Appendix~\ref{app:pseudocodes}): Start with a random dive through the tree to obtain $\theta^0, \zeta^0$, and $\kappa^0$ used to scale the features (see Table \ref{tab:state_features}).
\end{enumerate}

\begin{algorithm}[htbp!]
\caption{\textsc{$K$-B\&B-NodeSelection}}
\label{algorithm:k_bb_node_selection}
\DontPrintSemicolon
\KwInput{\quad Test instance $\P^{test}(N^{test})$, train instances $\P_1^{train}(N^{train}), \ldots, \P_I^{train}(N^{train})$ \\ \quad number of partitions for training $K^{train}$ testing $K^{test}$, level for training $L^{train}$ and testing \\ \quad $L^{test}$, quality threshold $\epsilon$, number of random dives per node $R$}
\KwOutput{\quad Objective value $\theta$, first-stage decisions $\x$, second-stage decisions $\y = \{\y_1, \ldots, \y_K \}$, \\
\quad subsets with scenarios $\barZ_k$ for all $k \in \{1,\ldots,K\}$}
\KwInit{\quad Incumbent partition: $\tau^i := \{\barZ_1, \ldots, \barZ_K\}$, where $\barZ_k = \emptyset$ for all $k \in \K$, \\
        \quad set containing all node partitions yet to explore: $\N := \{\tau^i\}$, \\
        \quad initial incumbent solution: $(\theta^i, \x^i, \y^i) := (\infty, \emptyset, \emptyset)$
}
\tcp{Preprocessing}
$\textit{model} \leftarrow \textsc{StrategyModel}(\P_1^{train}(N^{train}), \ldots, \P_I^{train}(N^{train}), K^{train}, L^{train}, \epsilon, R)$ \tcp*{Proc.~\ref{proc:strategy_model}}
$\textit{scaling info} \leftarrow \textsc{InitialRun}(\P^{test}(N^{test}), K^{test})$ \tcp*{Proc.~\ref{proc:init_run}}
\tcp{Tree search}
\While{$\N$ not empty}{
    \uIf{no solution yet or previous node pruned \label{alg_kbb_node:node_selection}}
        {Select a random node with partition $\tau = \{\barZ_1, \ldots, \barZ_K \}$ from $\mathcal{N}$, then $\mathcal{N} \leftarrow \mathcal{N}\setminus \{\tau\}$
        }
    \Else{
        $\tau := \{\barZ_1, \ldots, \barZ_{k^*} \cup \{\z^*\}, \ldots, \barZ_K\}$
    }    
    $(\theta^*, \x^*, \y^*) \leftarrow \textit{master problem}(\tau)$ \\
    \If{$\theta^* > \theta^i$}
    {
        Prune tree since current objective is worse than best solution found, continue to line \ref{alg_kbb_node:node_selection}.
 \label{alg_kbb_node:prune}\\ 
    }
    $(\z^*, \zeta^*) \leftarrow \textit{subproblem}(\theta^*, \x^*, \y^*)$ \\
    \uIf{$\zeta^* > 0$}
    {
        Solution not robust. Create $K$ new branches. \\ 
        \uIf{current level is more than $L^{test}$}
            {
            $k^* \leftarrow \textit{random uniform sample}([1, K])$ \\
            }
        \Else{
            Create feature vectors for $K$ child nodes. \\
            $(D_1, \ldots, D_K) \leftarrow \textit{generate features}(\textit{scaling info}, \theta, \zeta)$ \tcp*{steps of Section \ref{sec:features}}
            $k^* \leftarrow \textit{predict node qualities}(D_1, \ldots, D_K, \textit{model})$
            }
        Make $K$ new branches, of which the $k^*$-th is selected. \\
        \For{$k \in \{1, \ldots, K\} \setminus \{k^*\}$}
            {
                $\tau^k := \{\barZ_1, \ldots, \barZ_k \cup \{\z^*\}, \ldots, \barZ_K\}$ \\
                $\N \leftarrow \N \cup \{\tau^k\} $ \label{alg_kbb_node:branch}
            }
    }
    \Else{
        Current solution robust, prune tree. \label{alg_kbb_node:robust}\\
        $(\theta^i, \x^i, \y^i, \tau^i) \leftarrow (\theta^*, \x^*, \y^*, \tau)$
        }
}
\Return{$(\theta^i, \x^i, \y^i, \tau^i)$}
\end{algorithm}

\newpage
\section{Experiments} \label{sec:experiments}We now investigate if \textit{it is possible to learn a node selection strategy that generalizes (i) to other problem sizes, (ii) to different values of $K$, and (iii) to various problems.} We answer this question by a detailed study of two problems: capital budgeting (with loans) and shortest path \citep{subramanyam2020k}, whose formulations are given in Appendix~\ref{app:problems}. This section is set up as follows: First, the effectiveness of the original $K$-B\&B is tested on the problems. Then, we compare $K$-B\&B to \textsc{$K$-B\&B-NodeSelection}. We shall observe that the results obtained with our approach are very promising. 

For solving the MILPs, Gurobi 9.1.1 \citep{gurobi} is used. All computations of generating training data, $K$-B\&B, and \textsc{$K$-B\&B-NodeSelection} are performed on an Intel Xeon Gold 6130 CPU @ 2.1 GHz with 96 GB RAM. Training of the ML model is executed on an Intel Core i7-10610U CPU @ 1.8 GHz with 16 GB RAM. Our implementation along with the scripts to reproduce our results are available online\footnote{\url{https://github.com/estherjulien/KAdaptNS}}.

\subsection{Performance of $K$-B\&B}
In our experiments, we investigate the potential of improving the node selection strategy with ML. For that reason, it is important to identify problems on which such an improvement matters, \textit{i.e.}, where different partitions give varying outcomes and are nontrivial. To identify such problems, we run $K$-B\&B on several problems. Compared to the algorithm in \cite{subramanyam2020k}, we made some minor changes in $K$-B\&B:
\begin{itemize}
    \item Instead of breadth-first search, depth-first dives are performed with random node selection. {Depth-first search returns an incumbent solution early on, which is used for pruning nodes. This is increasingly vital when $K$ grows, as each node branches on $K$ other nodes. Likely a combination of breadth- and depth-first search would be preferred. This would require additional bookkeeping of the features for \textsc{$K$-B\&B-NodeSelection}.}
    \item In the starting node selection step (Step 2, Algorithm~\ref{algorithm:k_bb} in Appendix~\ref{app:pseudocodes}), a random node is taken from $\N$ instead of the first one. {This step is in line with random node selection.}
\end{itemize}
The quantity we are interested in is the relative change in the objective function value (OFV) -- the OFV of the first robust solution divided by the best one after 30 minutes. The higher the value, the more potential for a smart node selection strategy.

\begin{figure}[htbp!]
    \subfloat[Capital budgeting]{ \includegraphics[width=0.33\columnwidth]{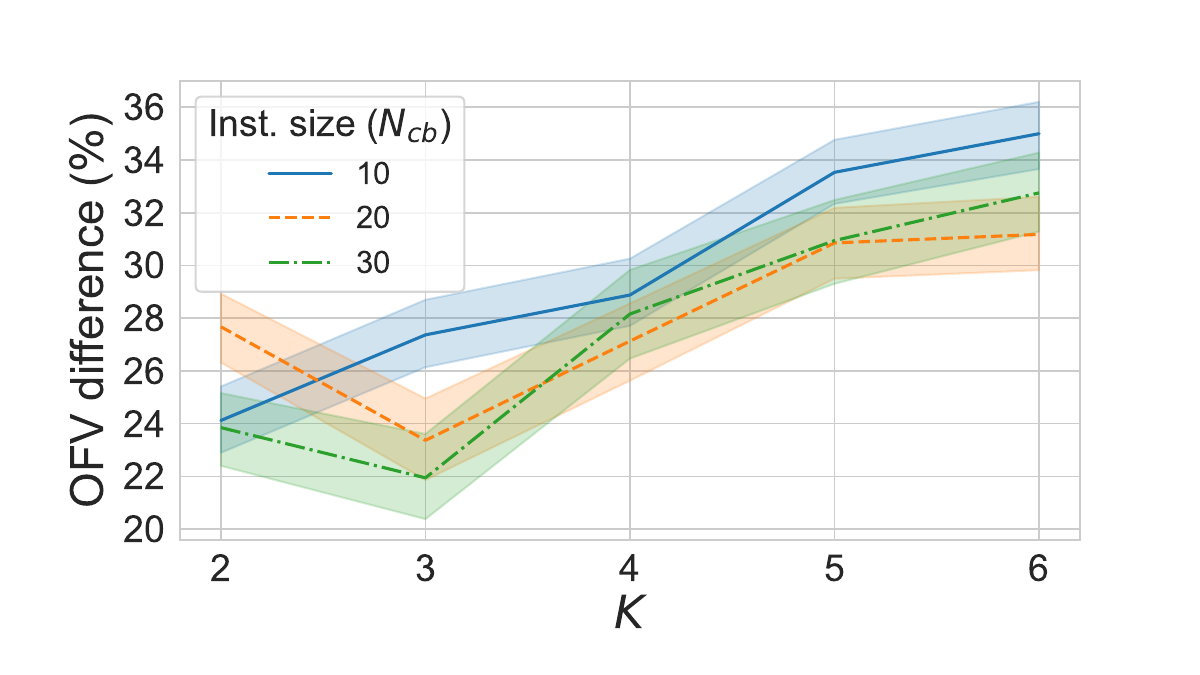} \label{fig:cb_diff}} 
    \subfloat[Shortest path (normal)]{\includegraphics[width=0.33\columnwidth]{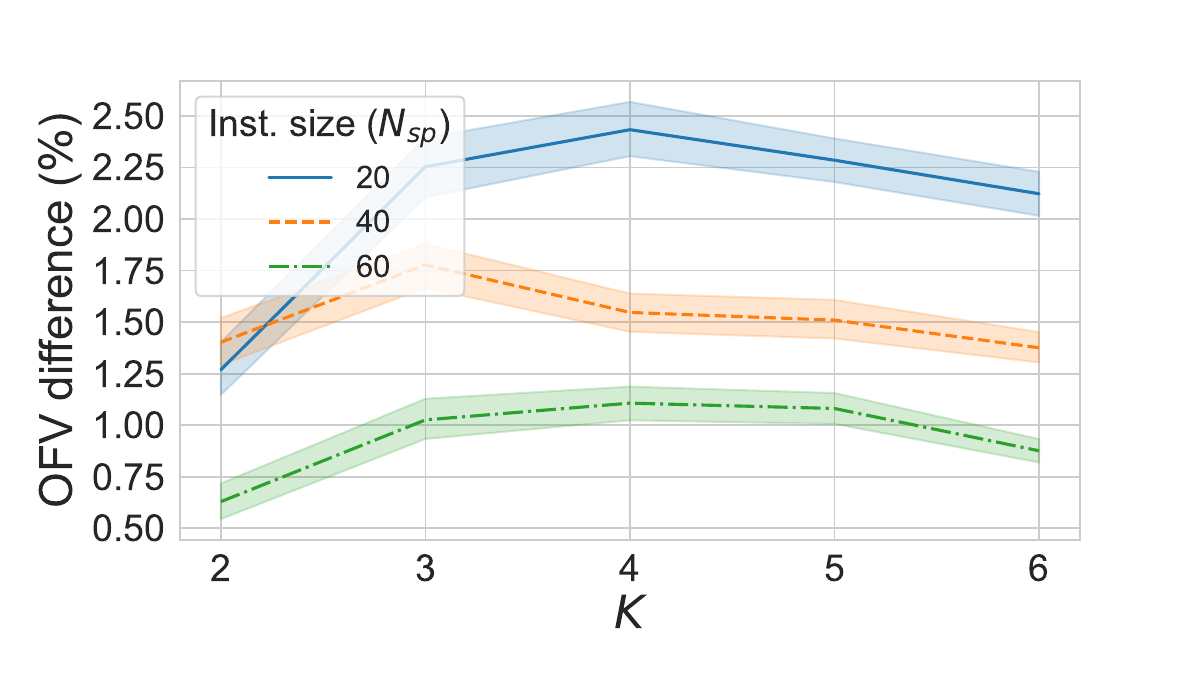} \label{fig:sp_normal_diff}}
    \subfloat[Shortest path (sphere)]{\includegraphics[width=0.33\columnwidth]{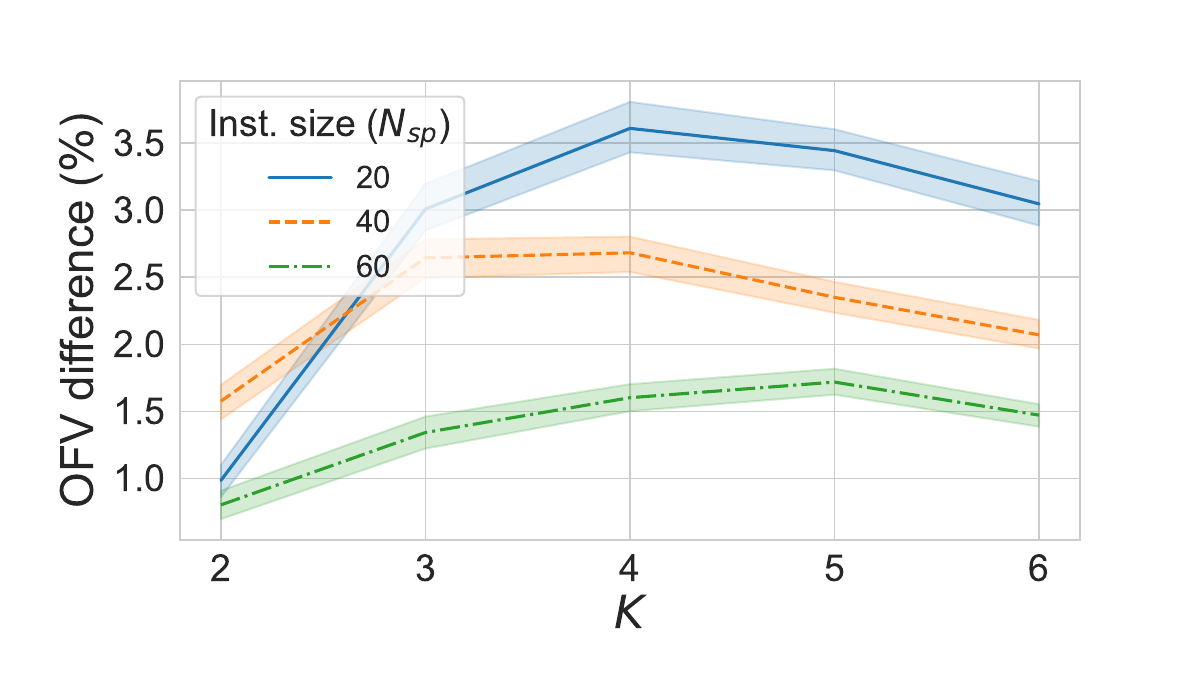} \label{fig:sp_sphere_diff}}
    \caption{Objective function value (OFV) difference (in $\%$) for capital budgeting and shortest path with `normal' and `sphere' instances, within 30 minutes using $K$-B\&B. For each problem type, experiments were done for 100 instances per $K \in \{2, \ldots, 6\}$ and instance size. The $75\%$ confidence interval (CI) is also given as shaded strips along the curves.} \label{fig:obj_diff_random}
\end{figure}

First, we consider the capital budgeting and the shortest path problems, only mentioned now and formally described later, for which the results are in Figure \ref{fig:obj_diff_random}. We observed that the objective function values of the capital budgeting instances are changing more than those of shortest path (in which nodes of a graph are located on a 2D plane). This is why we implemented another instance type for shortest path: graph with nodes located on a 3D sphere (see the Appendix). The OFV differences for these instances are still smaller than those of capital budgeting, but more than for the `normal' type. Therefore, further experiments on shortest path were conducted with the sphere instances alone.

Another problem on which we ran $K$-B\&B was the knapsack problem parameterized by the combination of the capacity ($c$), the number of items ($N_{ks}$), and the maximum deviation of the profit of items (\textit{i.e.}, the uncertainty parameter $\gamma$). For this problem, we considered instances of a similar size as the earlier problems, fixing $K = 4$ and $N_{ks} = 100$, and tested $K$-B\&B on 16 instances for all combinations of $c \in \{0.05, 0.15, 0.35, 0.5\}$ and $\gamma \in \{0.05, 0.15, 0.35, 0.5\}$. The OFV differences (in percentages) are given in Table \ref{tab:knapsack_perf}, where it is visible that different partitions barely play a role. Any random robust solution seems to be performing well. Thus, node selection will most likely not enhance $K$-B\&B for knapsack and will therefore not be tested.

\begin{table}[ht]
	\caption{The OFV difference (in $\%$) of the knapsack problem, within 30 minutes, for different values of the capacity ($c$) and uncertainty ($\gamma$) parameters. $K = 4$ and $N = 100$ are fixed.} 	\label{tab:knapsack_perf}
	\footnotesize
	\centering
	\begin{tabular}{l|llll}
        \toprule
         \multirow[c]{2}{*}{$c$} & \multicolumn{4}{c}{$\gamma$} \\ 
         {} &  0.05  &    0.15 &     0.35 &   0.5 \\
        \midrule
         0.05 &  0.002 &  0.086 &  0.000 &  0.000 \\
        0.15 &  0.002 &  0.011 &  0.040 &  0.000 \\
        0.35 &  0.001 &  0.001 &  0.019 &  0.038 \\
        0.5  &  0.001 &  0.003 &  0.012 &  0.027 \\
        \bottomrule
    \end{tabular}
\end{table}
Solving problem instances with $K$-B\&B takes a long time, where we have to think in the range of multiple hours until optimality is proven. This also holds for small-size instances. After having studied the convergence over runtime, we have fixed the time limit to 30 minutes for all the problem instances. For capital budgeting with $(N, K) = (10, 2)$ all instances could be solved within 30 minutes. For the same instance size and $K = 3$, $25$ out of $100$ are solved, only one for $K = 4$, and none for larger $K$. For the 3D shortest path problem with the smallest instance size and $K = 2$, only six instances are solved.

\subsection{Experimental setup \textsc{$K$-B\&B-NodeSelection}}
In our experiments, we also investigate what would be a good node selection strategy in $K$-B\&B that would perform well after it is trained on a selection of problems and then applied to other problems. Naturally, we are interested in generalization of trained tools to different instance sizes, $K$, and problem types. To this end, we have performed an ablation-study described in Table \ref{tab:experiment_types}, where each row informs about the similarity of the testing problem instances, compared to the training ones.

\begin{table}[htbp!]
	\caption{Different types of experiments (EXP), where the instance size $N$, $K$, and the problem itself can be different for training and testing.}	
	\label{tab:experiment_types}
	\footnotesize
	\centering
	\begin{tabular}{l|l|l|l}
		\toprule
		 \textbf{Type}    &   \textbf{Instance size $N$}   &    \textbf{$K$ }       &   \textbf{Problem}\\ 
        \midrule
         EXP1    &   Same                &   Same        &   Same \\
        EXP2    &   Same                &   Different   &   Same \\
        EXP3    &   Different           &   Same        &   Same \\
        EXP4    &   Different           &   Different   &   Same \\
        EXP5    &   Different           &   Same        &   Different \\
        EXP6    &   Different           &   Different   &   Different \\
		\bottomrule
	\end{tabular}
\end{table}

The problems we study are the capital budgeting problem with loans and the shortest path problem on a sphere (see the Appendix). We now describe the design choices we made regarding the ML model and data generation. As our focus lies in how ML is used for optimization and not in differences between ML models, we select a frequently-used model: random forest of \texttt{scikit-learn} \citep{scikit-learn} with default settings. We note that the training times do not exceed a couple of minutes for different data sets. 

As for the data generation process, it is governed by \textsc{StrategyModel} (Procedure~\ref{proc:strategy_model}), driven by the following parameters: $I$ (number of training instances), $L^{train}$ (depth level used in training data), and $R$ (number of dives per node). In these experiments, we instead made them depend on $T$ -- the total duration in hours, and $\iota$ -- the time per training instance in minutes. First, we set $I = 60T/\iota$. Selection of the right $L^{train}$ value is more challenging since for some problem instances the master problem takes a lot more time or deep trees are needed. Therefore, to get sufficiently many random dives for each node in $L^{train}$ within the time limit $\iota$, we make $L^{train}$ depend negatively on the total duration of the initial run (\textsc{InitialRun}, Procedure~\ref{proc:init_run}). This means that if a random dive takes a long time, the number of starting nodes should be lower, and therefore, the value of $L^{train}$ should also decrease. Finally, the number of dives per node ($R$) depends on how many nodes we have in level $L^{train}$ and the time we have for the instance ($\iota$).

First, the experiments of EXP1-EXP4 are run on the capital budgeting and shortest path problem. Then, the experiments where the training and testing problems are mixed, are conducted (EXP5 and EXP6). We only discuss a representative selection of the results in the main body, referring the reader to Appendix~\ref{app:results} for a complete overview. Finally, we discuss the feature importance scores we obtained with our random forests. 

\subsection{Capital budgeting}
\label{sec:experiments_cb}
The capital budgeting problem is a 2SRO problem where investments in a total of $N_{cb}$ projects can be made in two time periods. In the first period, the cost and revenue of these projects are uncertain. In the second period, these values are known but an extra penalty needs to be paid for postponement. The MILP formulation of capital budgeting has uncertain objective function, fixed number of uncertain constraints, and the dimension of $\mathcal{Z}$ is fixed to 4 for all instance sizes. For the full description, see Appendix~\ref{app:problems_cb}. Recall that \textsc{$K$-B\&B-NodeSelection} takes more parameters than the ones being tested in the ablation study: $L^{test}$ (level up to where node selection is performed), $\epsilon$ (node quality threshold), and $T$ and $\iota$ for generating training data. These parameters will be tuned in the first part of the experiments.

\subsubsection*{Parameter tuning (with EXP1)}
For this type of experiment, we use the same $K$ and $N$ for testing and training, applied to the smallest instance size: $N_{cb} = 10$, but for all $K \in \{2, 3, 4, 5, 6\}$. For different values of $K$, we noticed that different hyperparameters for generating training data were performing well. In Appendix~\ref{app:param_tuning_cb}, the tuning is performed of parameter $\iota$ (minutes spent per training instance) based on the number of data points obtained in total, together with some other information. The testing accuracy scores for different data sets we trained on, range between 0.92-0.99.

\begin{figure}[htbp!]
    \centering
    \includegraphics[width=.9\columnwidth]{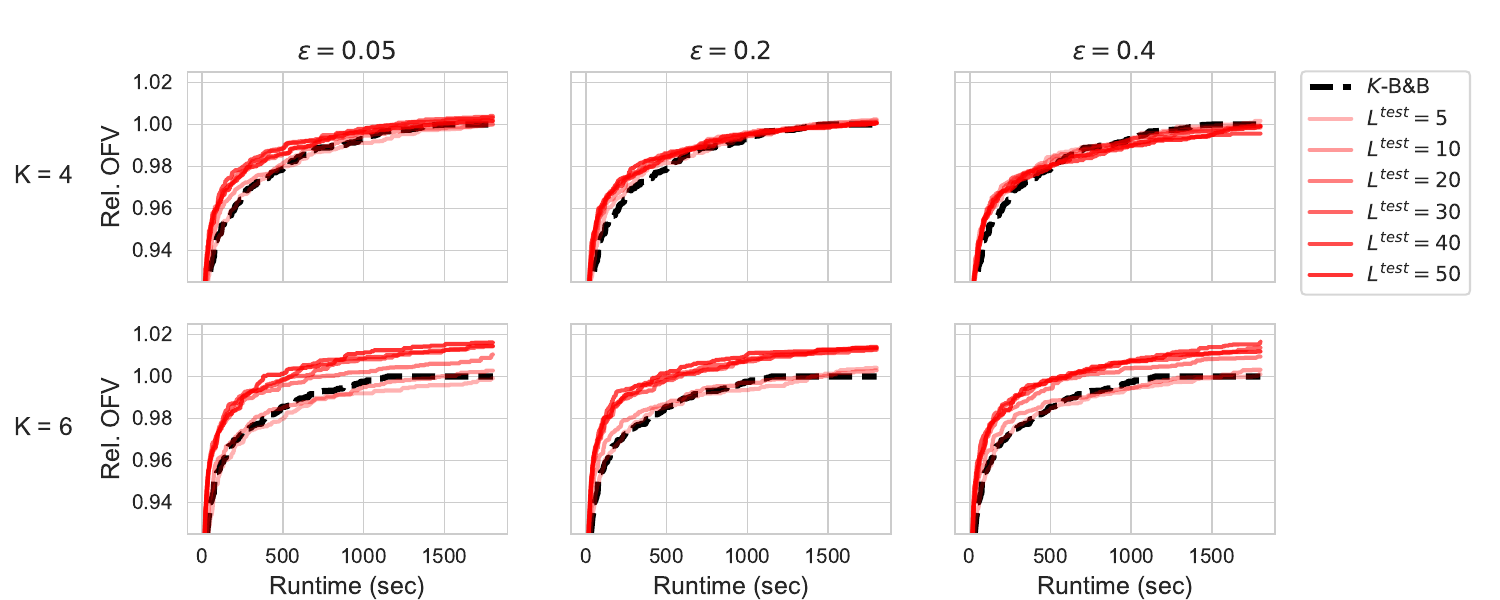}
    \caption{Results of $K$-B\&B with random dives and \textsc{$K$-B\&B-NodeSelection} with combinations of $K$, $L^{test}$ and $\epsilon$. The plots show the average relative objective value over the runtime of 100 instances for the capital budgeting problem, for $N = 10$.}
    \label{fig:cb_combinations_L_ct}
\end{figure}

We next tune the values of $L^{test}$ (level up to which node selection is performed) and $\epsilon$ (quality threshold) using the sets $L^{test} \in \{5, 10, 20, 30, 40, 50\}$ and $\epsilon \in \{0.05, 0.1, 0.2, 0.3, 0.4\}$.
Since there are 30 combinations per $K$ value, we only considered $K = 6$.
The results are shown in Figure \ref{fig:cb_combinations_L_ct}, where using a low value for $\epsilon$ and high $L^{test}$ gives the best results for both values of $K$. When a node $n$ has success probability $p_n$ larger than zero, this is already considered a good quality node. For further experiments of capital budgeting, we fix $\epsilon = 0.05$. Since high values of $L^{test}$ outperform lower ones, we also consider the possibility of applying the strategy always ($L^{test} = \infty$), with fixed $\epsilon = 0.05$. This option is analyzed in the Appendix. We noticed that choosing $L^{test} = \infty$ performs well for $K = 6$ but for lower values of $K$, choosing $L^{test} = 40$ gives even better solutions. Therefore, we continue with $L^{test} = 40$.

\begin{figure}[htbp!]
    \centering
        \includegraphics[width=\columnwidth]{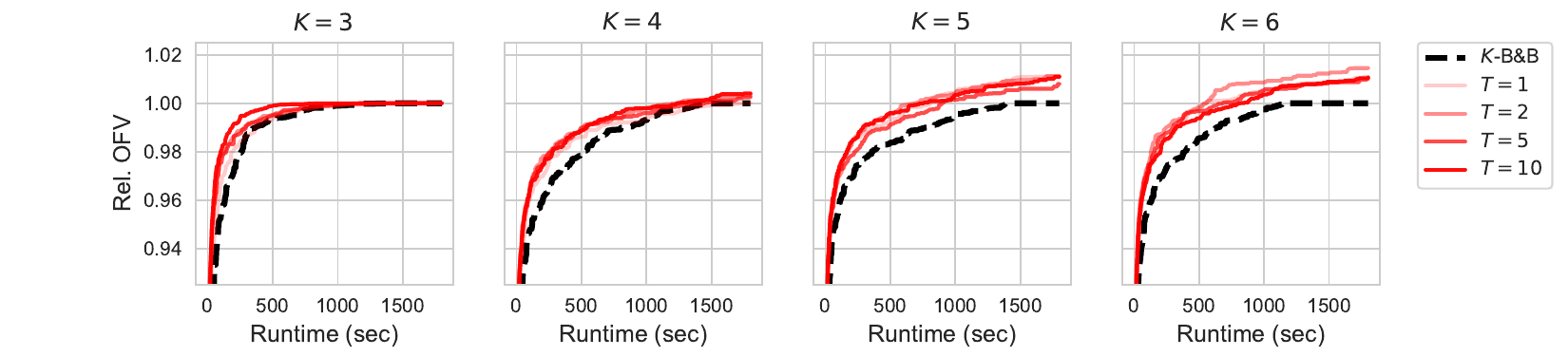}
    \caption{Results of $K$-B\&B with random dives and \textsc{$K$-B\&B-NodeSelection} with combinations of $K$ and $T$. The plots show the average objective over THE runtime of 100 instances for the capital budgeting problem, for $N = 10$.}
    \label{fig:cb_combinations_T}
\end{figure}

Another important parameter is the number of hours spent on generating training data ($T$). In Figure \ref{fig:cb_combinations_T}, results are shown for $T \in \{1, 2, 5, 10\}$ and $K\in \{3, 4, 5, 6\}$. Note that only for $K = 3$ higher values of $T$ result in a better performance of \textsc{$K$-B\&B-NodeSelection}. For the other values of $K$ the performance is similar, or even worse, for higher values of $T$ compared to lower ones. For further experiments, we shall use $T = 2$.

{
\subsubsection*{Alternative scaling}
The step taken to obtain scales for the objective ($\theta^0$), violation ($\zeta^0$), and depth ($\kappa^0$) in \textsc{InitialRun} seems to be avoidable if we decide to take other values for these parameters. In this section, we will describe a small experimental study of different scaling methods: (i) \textsc{InitialRun}, (ii) Alternative, and (iii) No scaling, where the latter corresponds to setting $\theta^0=\zeta^0=\kappa^0=1$. The Alternative setting takes as scale for $\theta^0$ the objective of the root node, and $\kappa^0=K^{\log(\dim(\mathcal{Z}))}$. Next to these `internal' scaling approaches, one can also use min/max scaling of all feature values. The scales are taken from the training data set and then used in the testing phase. These two scaling approaches can be used in combination. The results obtained for the capital budgeting problem, with $N=10$ and $K^{test}=K^{train}=6$, are given in Table \ref{tab:scaling_test}. Here, for the first time, we use performance statistics to measure the effectiveness of the method. From the results presented in the table, we can conclude that, on average, solutions found after 30 minutes are best when \textsc{InitialRun} is used. Additional min/max scaling on top of the different scaling approaches gives mixed results.}
        \begin{table}[htbp!]
            \centering
            \footnotesize
            \caption{The results for a variety of scaling methods tested on capital budgeting problem, with $N=10$, $K^{test}=K^{train}=6$, on 16 instances. Different `internal' scaling methods: `\textsc{InitialRun}', `Alternative', and `No scaling'. These internal scaling methods can be combined with min/max scaling. The results for all combinations is given. Four performance statistics are given: (i) `OFV 30m': Difference (in \%) in relative OFV found after 30 minutes from \textsc{$K$-B\&B} to \textsc{$K$-B\&B-NodeSelection} (higher is better). (ii) `OFV 1m': is the relative OFV from the two methods after one minute. (iii): `to OFV=1': Runtime speedup (in \%) to reach a relative OFV of 1 for \textsc{$K$-B\&B-NodeSelection} compared to \textsc{$K$-B\&B} (higher is better). (iv) `NS to OFV=1': Number of instances of \textsc{$K$-B\&B-NodeSelection} that reached a relative OFV of 1 (higher is better).
            }
            \label{tab:scaling_test}
            \begin{tabular}{lr|rr|rr|rr}
            \toprule
              \multirow[c]{2}{*}{\textbf{Statistic}} & & \multicolumn{2}{r|}{\textbf{Initial dive}} & \multicolumn{2}{r|}{\textbf{Alternative}} & \multicolumn{2}{r}{\textbf{No scaling}} \\ 
            & \textbf{min/max} & NO & YES & NO & YES & NO & YES \\
            \midrule
            OFV 30m & & 2.10 & 1.98 & 1.28 & 1.11 & 1.33 & 1.50 \\
            OFV 1m & & 2.14 & 2.39 & 1.49 & 1.93 & 2.03 & 1.39 \\
            to OFV=1 & & 44.10 & 42.53 & 52.58 & 38.08 & 51.72 & 44.29 \\
             NS to OFV=1 & & 16 & 16 & 13 & 14 & 13 & 14 \\
            \bottomrule
            \end{tabular}
        \end{table}

\subsubsection*{Results}
{We compare $K$-B\&B to \textsc{$K$-B\&B-NodeSelection} for all combinations of $K^{test}, K^{train} \in \{2, 3, 4, 5, 6\}$, and $N_{test} \in \{10, 20, 30\}$. The results of the full set of combinations are displayed in Appendix~\ref{app:results_cb}. Table~\ref{tab:cb_exp1_exp2} lists the results using four performance statistics of EXP1 and EXP2 (described in the caption). The table shows that \textsc{$K$-B\&B-NodeSelection} outperforms \textsc{$K$-B\&B}, with speedups ranging from on average 2 to 50\%, except for $K^{train}=2$ and the instances of $K^{test}=2$. Higher values of $K^{train}$ generally perform better both in terms of relative OFV and speedup.
\\
In the next sections, we look at specific problem instances in more detail, where we also focus on the stability of the algorithm. Due to the many combinations for EXP3-EXP4, a table would be too convoluted.}

\begin{table}[htbp!]
    \centering
    \footnotesize
    \caption{Combined results of EXP1 (along diagonal) and EXP2 for the capital budgeting problem. The four statistics described in the caption of Table \ref{tab:scaling_test} are given: (i) `OFV 30m' (higher is better), (ii) `OFV 1m', (iii) `to OFV=1' (higher is better), and (iv) `NS to OFV=1' (higher is better).}
    \label{tab:cb_exp1_exp2}
\begin{tabular}{c|l|rrrrr}
\toprule
 \multirow[c]{2}{*}{\textbf{$K^{test}$}} & \multirow[c]{2}{*}{\textbf{Statistic}} & \multicolumn{5}{c}{\textbf{$K^{train}$}}\\  
 &  & 2 & 3 & 4 & 5 & 6 \\
 \midrule
\multirow[c]{4}{*}{2} & OFV 30m & -0.00 & -0.00 & -0.00 & -0.00 & -0.00 \\
 & OFV 1m & -0.13 & -0.12 & -0.09 & 0.02 & -0.09 \\
 & to OFV=1 & -101.47 & -78.31 & -60.72 & -54.80 & -33.46 \\
 & NS to OFV=1 & 47 & 46 & 56 & 50 & 55 \\
\midrule
\multirow[c]{4}{*}{3} & OFV 30m & -0.08 & -0.13 & 0.05 & -0.06 & 0.12 \\
 & OFV 1m & 0.56 & 2.07 & 2.38 & 2.31 & 2.74 \\
 & to OFV=1 & 3.13 & 18.67 & 26.95 & 12.00 & 24.38 \\
 & NS to OFV=1 & 47 & 44 & 49 & 56 & 64 \\
\midrule
 \multirow[c]{4}{*}{4} & OFV 30m & -0.47 & -0.02 & 0.34 & 0.66 & 0.95 \\
 & OFV 1m & -0.19 & 1.58 & 1.98 & 2.22 & 2.33 \\
 & to OFV=1 & 21.19 & 35.19 & 43.33 & 45.80 & 49.88 \\
 & NS to OFV=1 & 33 & 45 & 62 & 74 & 78 \\
\midrule 
 \multirow[c]{4}{*}{5} & OFV 30m & -0.08 & 0.43 & 0.66 & 1.15 & 1.55 \\
 & OFV 1m & -0.70 & 1.36 & 1.34 & 1.38 & 1.79 \\
 & to OFV=1 & -0.57 & 35.43 & 40.61 & 48.15 & 52.80 \\
 & NS to OFV=1 & 53 & 66 & 69 & 79 & 88 \\
\midrule 
 \multirow[c]{4}{*}{6} & OFV 30m & -0.02 & 0.74 & 0.88 & 1.18 & 1.74 \\
 & OFV 1m & 0.07 & 1.49 & 1.68 & 1.68 & 2.33 \\
 & to OFV=1 & 13.50 & 30.99 & 24.11 & 27.75 & 31.75 \\
 & NS to OFV=1 & 49 & 69 & 76 & 78 & 89 \\
 \bottomrule
\end{tabular}
\end{table}

\begin{figure}[htbp!]
    \centering
    \includegraphics[width=.9\columnwidth]{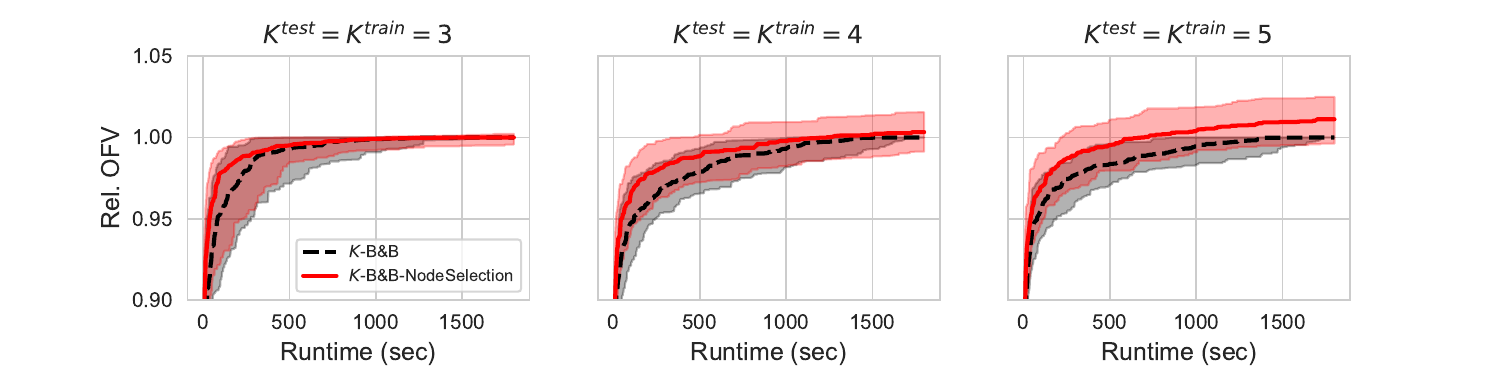} 
    \caption{EXP1 results for $K^{train} = K^{test} \in \{3, 4, 5\}$  on 100 instances. {The black line gives the average relative objective function value (Rel. OFV) over the runtime (in seconds) of \textsc{$K$-B\&B}, with a 30 minute time limit. The red line is the Rel. OFV trajectory of \textsc{$K$-B\&B-NodeSelection}. The shaded area around the lines is their respective 75\% confidence interval.}}
    \label{fig:cb_exp1}
\end{figure}

\begin{figure}[htbp!]
    \centering
    \includegraphics[width=.9\columnwidth]{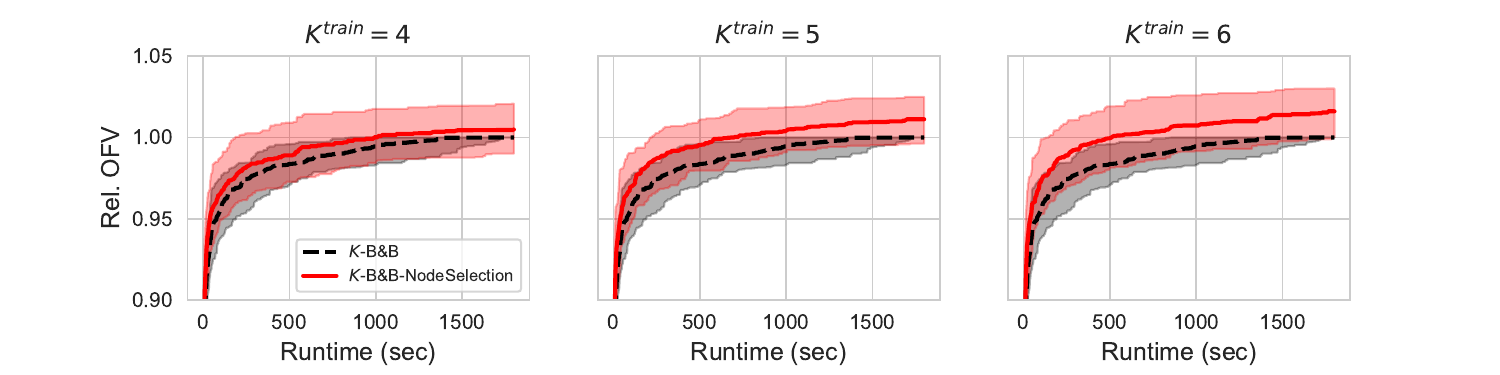} 
    \caption{EXP2 results for $K^{test}=5$ and $K^{train} \in \{4, 5, 6\}$ on 100 instances. {The black line gives the average Rel. OFV of \textsc{$K$-B\&B} and the red line that of \textsc{$K$-B\&B-NodeSelection}.}}
    \label{fig:cb_exp2}
\end{figure}

\paragraph{EXP1 and EXP2 results.}
For EXP1, the results for $K^{train} = K^{test} \in \{3, 4, 5\}$ are shown in Figure \ref{fig:cb_exp1}. As in Figure \ref{fig:obj_diff_random}, the shaded strips around the solid curves show the confidence interval. We can see that \textsc{$K$-B\&B-NodeSelection} outperforms $K$-B\&B for all $K$. Moreover, the convergence is steeper: good solutions are found earlier. For $K \in \{4, 5\}$ the final solution is also better when node selection is guided by ML predictions. Then for EXP2, where we also apply ML models that are trained with $K^{train} \neq K^{test}$, we see the performance of \textsc{$K$-B\&B-NodeSelection} improving when $K^{train}$ increases. This indicates that data obtained with higher values of $K^{train}$ are more informative than those with lower values. See Figure \ref{fig:cb_exp2} for an illustration with $K^{test} = 5$ and $K^{train} \in \{4, 5, 6\}$. 

\begin{figure}[htbp!]
    \centering
        \includegraphics[width=\columnwidth]{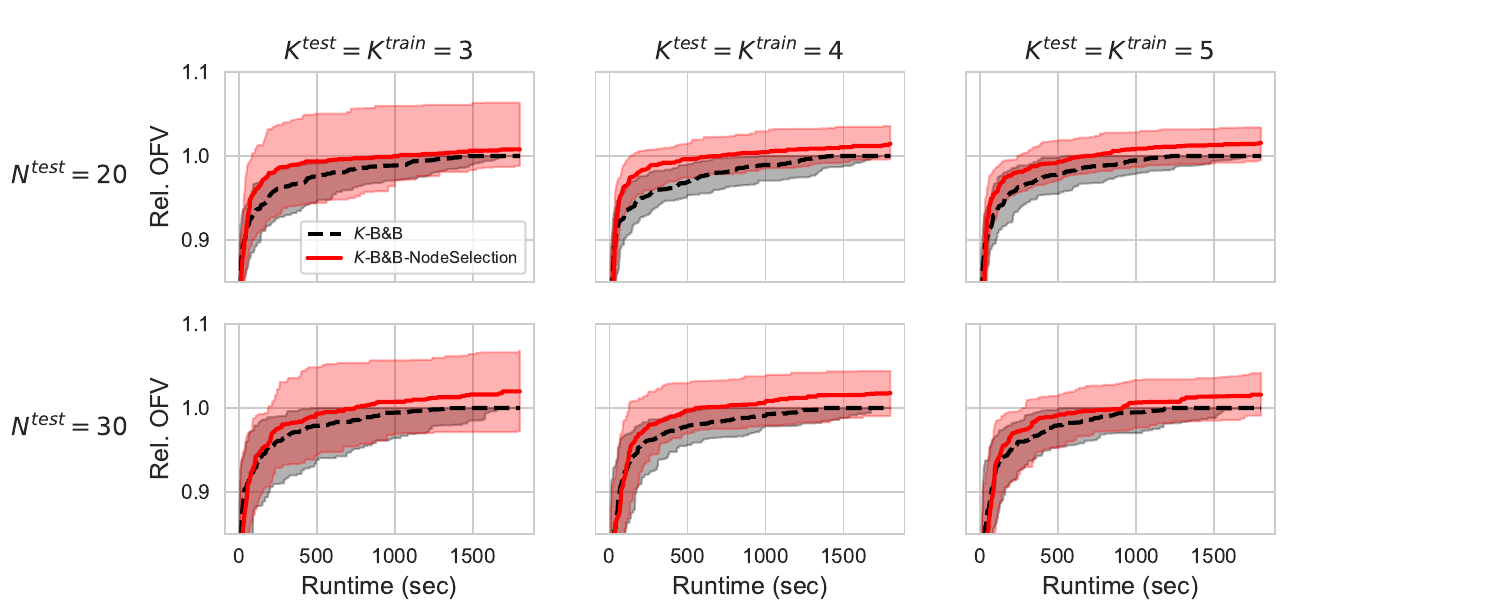}
    \caption{EXP3 results for $K^{train} = K^{test} \in \{3, 4, 5\}$ and $N^{test} \in \{20, 30\}$ on 100 instances. {The black line gives the average Rel. OFV of \textsc{$K$-B\&B} and the red line that of \textsc{$K$-B\&B-NodeSelection}.} }
    \label{fig:cb_exp3}
\end{figure}

\begin{figure}[htbp!]
    \centering
        \includegraphics[width=\columnwidth]{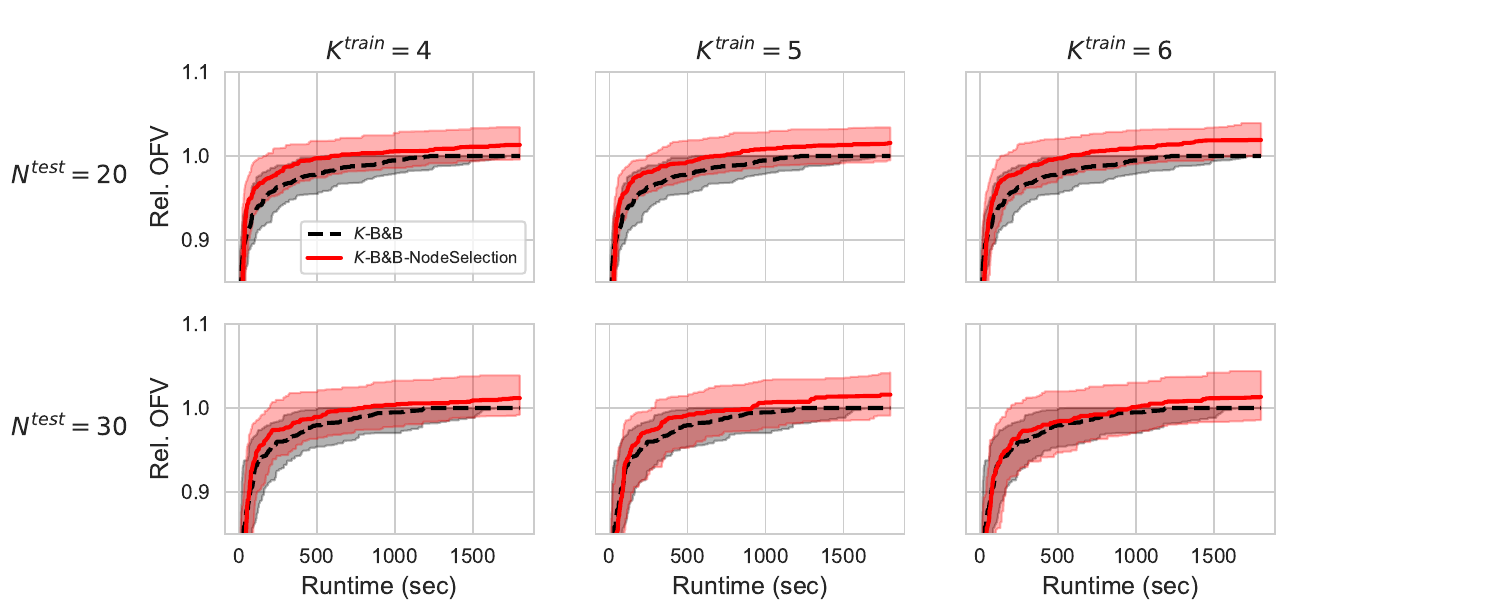}
    \caption{EXP4 results for $K^{test} = 5$, $K^{train}  \in \{4, 5, 6\}$, and $N^{test} \in \{20, 30\}$ on 100 instances. {The black line gives the average Rel. OFV of \textsc{$K$-B\&B} and the red line that of \textsc{$K$-B\&B-NodeSelection}.}}
    \label{fig:cb_exp4}
\end{figure}

\paragraph{EXP3 and EXP4 results.} In Figure \ref{fig:cb_exp3}, we depict the results for experiments with same $K$ but different instance sizes. {On average, \textsc{$K$-B\&B-NodeSelection} performs better both in terms of speedup and final relative OFV found: for $N^{test}=20$, the speedups are between 20-50\% and relative OFV increase around 1.5\%. For $N^{test}=30$, we see speedups ranging from 11-47\% and relative OFV increases of 1.6-2.1\%.} However, the confidence interval of $N^{test}=30$ is larger than that of $N^{test} = 10$. This suggests that testing on higher values of $N$ gives rise to a higher risk. Testing on higher instance sizes for different values of $K$ (Figure \ref{fig:cb_exp4}) has a similar effect as before: higher values of $K^{train}$ result in better performance {(Rel. OFV ranging from 1.4-1.9\%, and speedup 5-27\%)}, although marginally for some parameter combinations.

\subsection{Shortest path on a sphere} \label{sec:experiments_sp}
The shortest path problem can be described as a 2SRO problem where we make the planning of a route from source to target, for which the lengths of the $N_{sp}$ arcs are uncertain. This problem only has second-stage decisions and an uncertain objective. The dimension of the uncertainty set grows with the number of arcs in the graph. For the full description, see Appendix~\ref{app:problems_sp}. Since there is only a second stage, some attributes disappear for this problem: `Deterministic first-stage' (attribute 6), and the static objective problem-related attributes (8 and 9). Hence, we are left with six attributes for this problem.

\subsubsection*{Parameter tuning (with EXP1)} 
\begin{figure}[htbp!]
    \centering
        \includegraphics[width=\columnwidth]{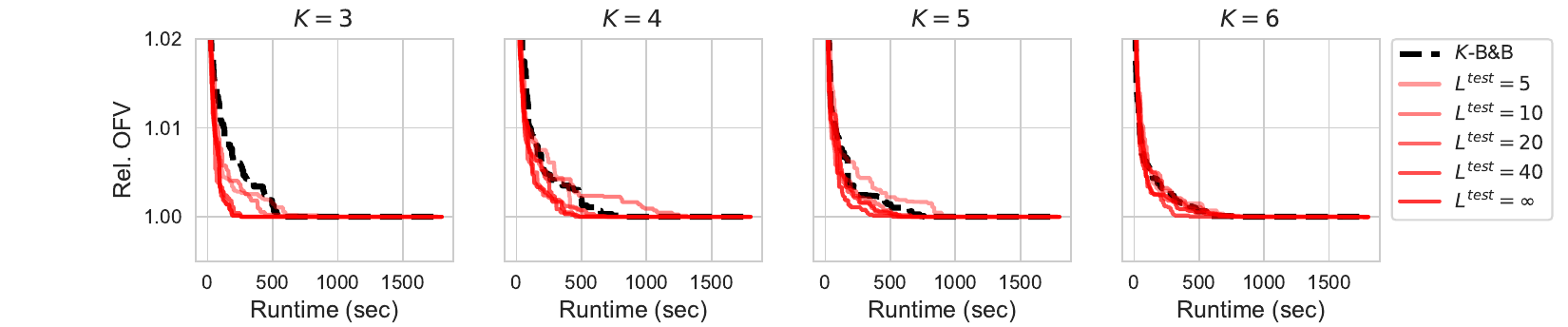} \caption{Results of $K$-B\&B with random dives and \textsc{$K$-B\&B-NodeSelection} with combinations of $K$ and $L^{test}$. The plots show the average objective over runtime of 100 instances for the shortest path problem, for $N = 20$. }
    \label{fig:sp_diffL}
\end{figure}

As we did for the capital budgeting problem, we first tune the parameter $\iota$. The details of this parameter tuning are given in Appendix~\ref{app:param_tuning_sp}. The testing accuracy scores for different data sets we trained on is lower than for capital budgeting; between 0.88-0.97. Since the distribution of the probability success values $p_n$ is very similar to the capital budgeting problem, the quality threshold is set to $\epsilon = 0.05$. In Figure \ref{fig:sp_diffL}, the level is tested again with $ L^{test} \in \{5, 10, 20, 40, \infty\}$. Here we see that especially for smaller values of $K$, a higher level of $L^{test}$ outperforms others. Therefore, we select $L^{test} = \infty$ for the remainder of the experiments.

\begin{figure}[htbp!]
    \centering
        \includegraphics[width=\columnwidth]{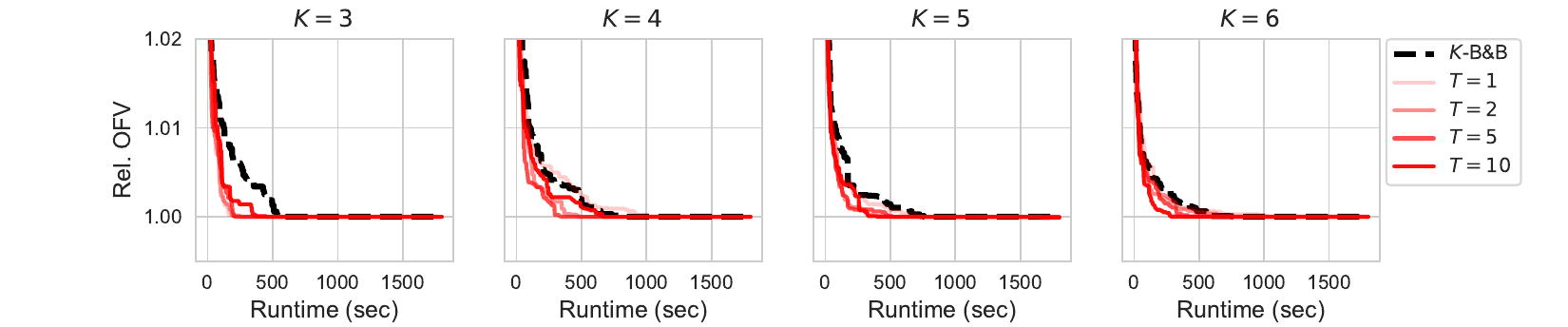} \caption{Results of $K$-B\&B with random dives and \textsc{$K$-B\&B-NodeSelection} with combinations of $K$ and $T$. The plots show the average objective over a runtime of 100 instances for the shortest path problem, for $N = 20$. }
    \label{fig:sp_diffT}
\end{figure}

Also note that when $K$ grows, the performances of $K$-B\&B and \textsc{$K$-B\&B-NodeSelection} are very similar. For the shortest path problem, we noticed that more training data points led to a substantial performance gain. See Figure \ref{fig:sp_diffT} for these results. Therefore, for EXP1-EXP4, we select $T$ for each $K$ separately. For an overview of which parameters are chosen per $K$; see Appendix~\ref{app:param_tuning_sp}.

{\subsubsection*{Results}
The entirety of the results for all combinations of $K^{test}$, $K^{train}$, and $N^{test}$ for $K$-B\&B versus \textsc{$K$-B\&B-NodeSelection} can be found in Appendix~\ref{app:results_sp}. Again, some performance statistics for EXP1-EXP2 are given, see Table \ref{tab:sp_exp1_exp2}. This table shows that the relative OFV found by \textsc{$K$-B\&B-NodeSelection} is marginally better than those of \textsc{$K$-B\&B} (0.01-0.12\% when excluding $K^{train}=2$), but that the speedups can be very significant: for $K = 6$, the average speedups range from 53-94\%. For this problem class, the instances of $K=2$ seem to benefit from ML enhancements.}

\begin{table}[htbp!]
    \centering
    \footnotesize
    \caption{Combined results of EXP1 (along diagonal) and EXP2 for the shortest path problem. The four statistics described in the caption of Table \ref{tab:scaling_test} are given: (i) `OFV 30m' (higher is better), (ii) `OFV 1m', (iii) `to OFV=1' (higher is better), and (iv) `NS to OFV=1' (higher is better).}
    \label{tab:sp_exp1_exp2}
\begin{tabular}{c|l|rrrrr}
\toprule
 \multirow[c]{2}{*}{\textbf{$K^{test}$}} & \multirow[c]{2}{*}{\textbf{Statistic}}  & \multicolumn{5}{c}{\textbf{$K^{train}$}}\\  
 &  & 2 & 3 & 4 & 5 & 6 \\
 \midrule
\multirow[c]{4}{*}{2} & OFV 30m & 0.05 & 0.07 & 0.03 & 0.03 & 0.08 \\
 & OFV 1m & 0.08 & 1.10 & 1.08 & 1.04 & 0.10 \\
 & to OFV=1 & 36.63 & 15.83 & 24.77 & 22.31 & 12.31 \\
 & NS to OFV=1 & 98 & 100 & 99 & 100 & 97 \\
\midrule
 \multirow[c]{4}{*}{3} & OFV 30m & 0.01 & 0.14 & 0.13 & 0.11 & 0.12 \\
 & OFV 1m & 0.17 & 0.45 & 0.51 & 0.19 & 0.52 \\
 & to OFV=1 & 60.44 & 46.93 & 34.30 & 49.02 & 29.06 \\
 & NS to OFV=1 & 99 & 100 & 99 & 100 & 100 \\
\midrule
 \multirow[c]{4}{*}{4} & OFV 30m & -0.12 & 0.10 & 0.12 & 0.11 & 0.12 \\
 & OFV 1m & 0.89 & 0.41 & 0.41 & 0.47 & 0.65 \\
 & to OFV=1 & 55.90 & 63.34 & 75.11 & 45.12 & 45.79 \\
 & NS to OFV=1 & 99 & 100 & 100 & 100 & 100 \\
\midrule
 \multirow[c]{4}{*}{5} & OFV 30m & -0.07 & 0.12 & 0.09 & 0.15 & 0.10 \\
 & OFV 1m & -0.33 & 0.21 & 0.16 & 0.23 & 0.37 \\
 & to OFV=1 & 62.87 & 81.84 & 77.14 & 64.61 & 52.25 \\
 & NS to OFV=1 & 100 & 100 & 100.00 & 100 & 100 \\
\midrule 
 \multirow[c]{4}{*}{6} & OFV 30m & -0.07 & 0.04 & 0.01 & 0.10 & 0.07 \\
 & OFV 1m & -0.07 & 0.11 & 0.17 & 0.20 & 0.23 \\
 & to OFV=1 & 52.26 & 67.37 & 67.55 & 88.17 & 93.68 \\
 & NS to OFV=1 & 100 & 100 & 100 & 100 & 100 \\
 \bottomrule
\end{tabular}
\end{table}

\begin{figure}[htbp!]
    \centering
    \includegraphics[width=.9\columnwidth]{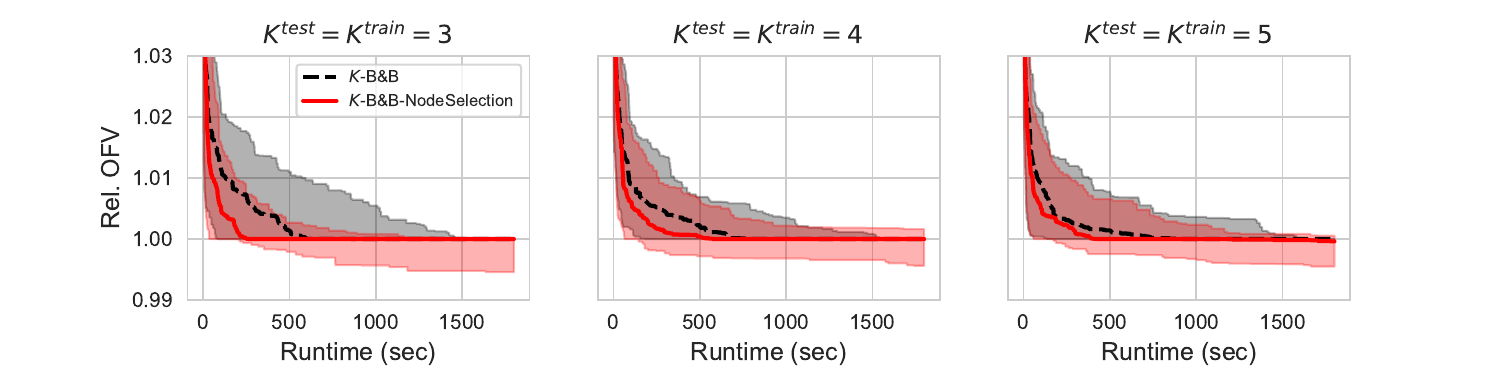}
    \caption{EXP1 results for $K^{train} = K^{test} \in \{3, 4, 5\}$ on 100 instances. {The black line gives the average relative objective function value (Rel. OFV) over the runtime (in seconds) of \textsc{$K$-B\&B}, with a 30 minute time limit. The red line is the Rel. OFV trajectory of \textsc{$K$-B\&B-NodeSelection}. The shaded area around the lines is their respective 80\% confidence interval.}}
    \label{fig:sp_exp1}
\end{figure}

\begin{figure}[htbp!]
    \centering
    \includegraphics[width=.9\columnwidth]{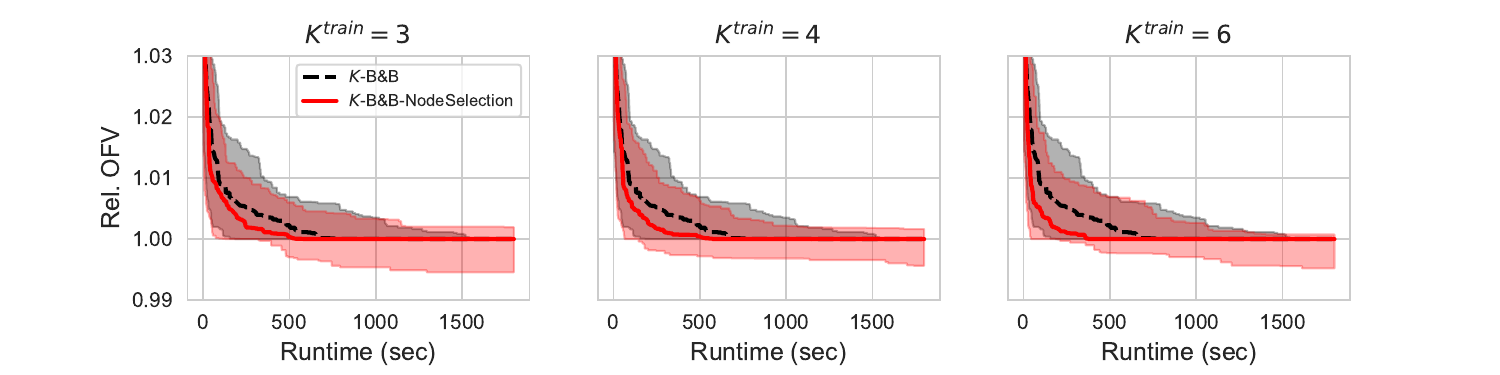}
    \caption{EXP2 results for $K^{test} = 4$ and $K^{train} \in \{3, 4, 6\}$ on 100 instances. {The black line gives the average Rel. OFV of \textsc{$K$-B\&B} and the red line that of \textsc{$K$-B\&B-NodeSelection}.}}
    \label{fig:sp_exp2}
\end{figure}

\paragraph{EXP1 and EXP2 results.} As is visible from Figure \ref{fig:sp_exp1}, even though no better solutions are found (probably because the optimal solution is found early on), it is still noticeable that \textsc{$K$-B\&B-NodeSelection} converges faster than $K$-B\&B. This phenomenon is thus consistent over the two problems. The convergence of $K$-B\&B over different values of $K^{test}$ differs significantly (\textit{e.g.}, compare $K^{test}$ equal to 3 and 5). The convergence of \textsc{$K$-B\&B-NodeSelection} is however quite stable across different values of $K^{test}$. Then, for EXP2, in Figure \ref{fig:sp_exp2}, we see that for $K^{test} = 4$, $K^{train}=6$ performs best. For other values of $K^{test}$, $K^{train}=6$ behaves well too, just as it did for the capital budgeting problem.

\begin{figure}[htbp!]
    \centering
        \includegraphics[width=\columnwidth]{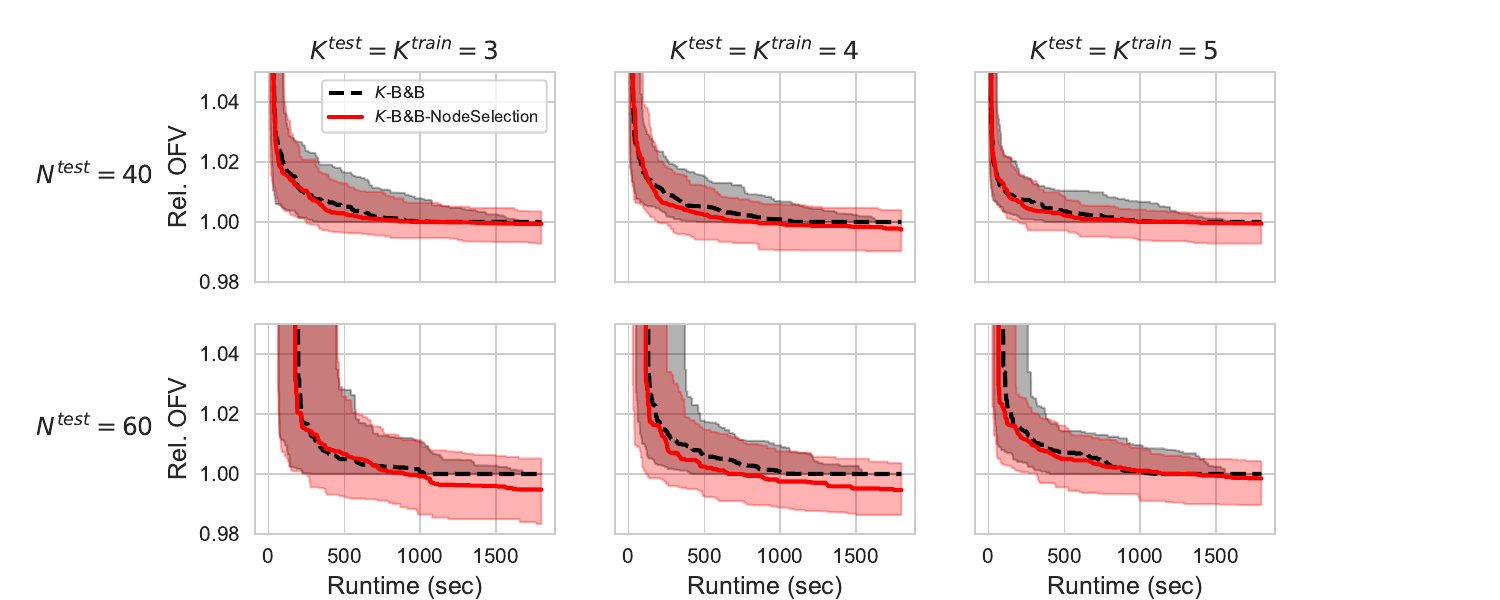}
    \caption{EXP3 results for $K^{train} = K^{test} \in \{3, 4, 5\}$ and $N^{test} \in \{40, 60\}$ on 100 instances. {The black line gives the average Rel. OFV of \textsc{$K$-B\&B} and the red line that of \textsc{$K$-B\&B-NodeSelection}.}}
    \label{fig:sp_exp3}
\end{figure}

\begin{figure}[htbp!]
    \centering
        \includegraphics[width=\columnwidth]{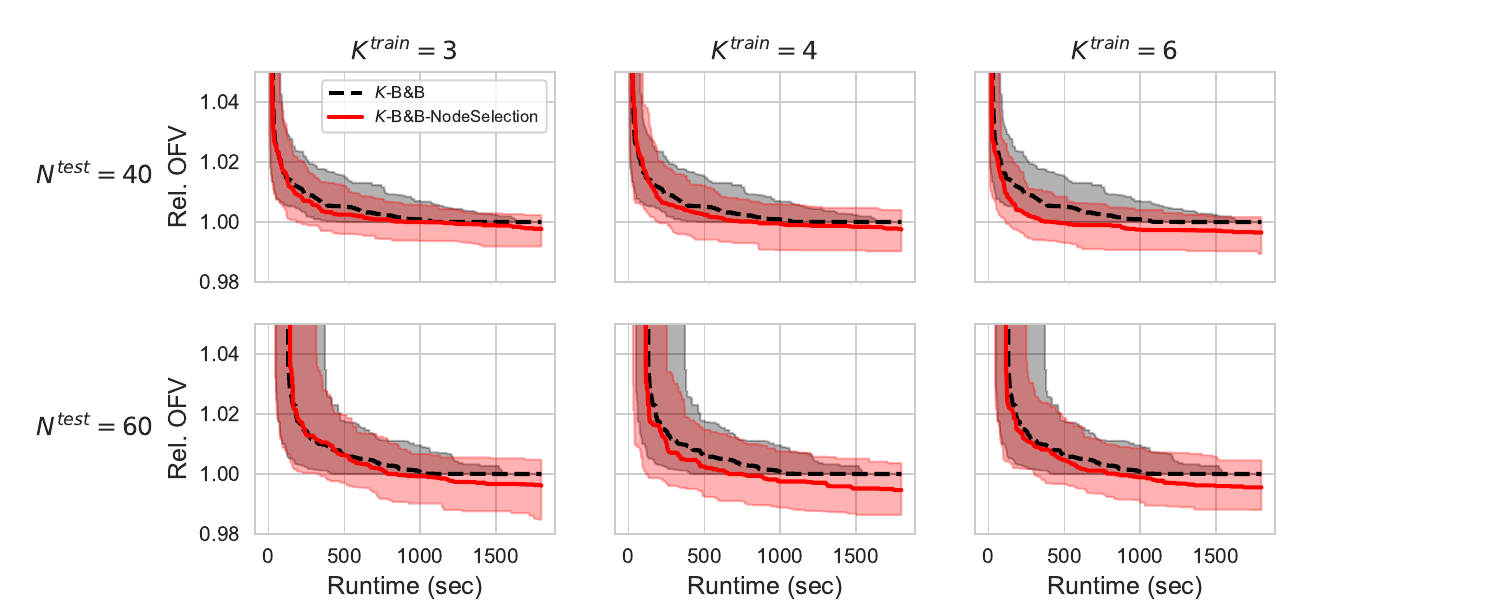}
    \caption{EXP4 results for $K^{test} = 4$, $K^{train} \in \{3, 4, 6\}$, and $N^{test} \in \{40, 60\}$ on 100 instances. {The black line gives the average Rel. OFV of \textsc{$K$-B\&B} and the red line that of \textsc{$K$-B\&B-NodeSelection}.}}
    \label{fig:sp_exp4}
\end{figure}

\paragraph{EXP3 and EXP4 results.}  We see in Figure \ref{fig:sp_exp3} that the two algorithms behave very similarly. {The average relative OFV ranges from 0.14-0.55\%. Interestingly, as not easily visible from the figure, the average speedup values are quite high: ranging from 21-98\%. The speedup of 98\% is achieved for $K = 6$ and $N=40$, where the line of \textsc{$K$-B\&B} consistently lies slightly above the one of \textsc{$K$-B\&B-NodeSelection}.} In terms of stability, the confidence interval of \textsc{$K$-B\&B-NodeSelection} grows with $N^{test}$. However, this does not necessarily mean that testing and training on different sizes is not stable: the CI of \textsc{$K$-B\&B} is also bigger. Moreover, the CI is mostly in the region below one, which indicates that we mainly have well-performing outliers. In Figure \ref{fig:sp_exp4}, the results of EXP4 are shown, {where the average relative OFV improvement ranges from 0.29-0.51\% and the speedups ($K=4$ excluded) are between 32-77\%}. We see that for $N^{test}=40$, $K^{train}=6$ performs best, but not necessarily for the biggest instance size.

\subsection{Training and testing on different problems}
\label{sec:experiments_mix}
In this section, we show the results of EXP5 and EXP6, where we apply the node selection strategy to a different problem than it has been trained on. Note that the shortest path problem does not have first-stage decisions. This results in the features being slightly different than for the capital budgeting problem. Therefore, to create a model that can be trained by shortest path data, and used for the capital budgeting problem, the first-stage-related attributes are not constructed while running \textsc{$K$-B\&B-NodeSelection}. Full results of these experiments are given in Appendix~\ref{app:results_mixed}. 

Figure \ref{fig:cb_exp5} shows the results of EXP5 for the capital budgeting problem. We can see that the performances of the two algorithms are very close. Recall that $K^{train} = 6$ previously resulted in (one of) the best solutions. Then, when we look at some of the solutions of EXP6 with $K^{train} = 6$ (see Figure \ref{fig:cb_exp6}), we see this is not the case here.

\begin{figure}[htbp!]
    \centering
    \subfloat[EXP5]{
        \includegraphics[width=.9\columnwidth]{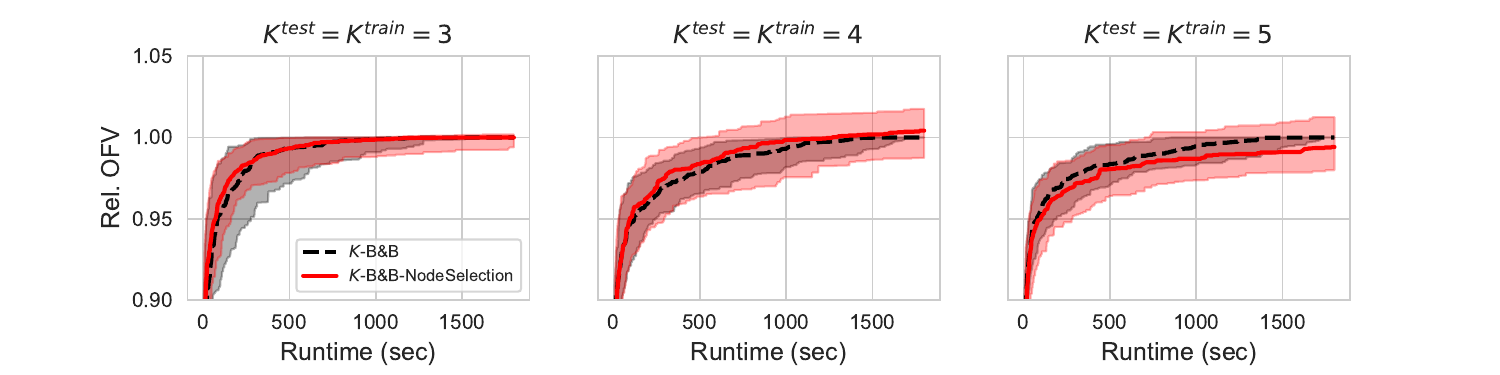} \label{fig:cb_exp5}} \\ 
    \subfloat[EXP6]{        \includegraphics[width=.9\columnwidth]{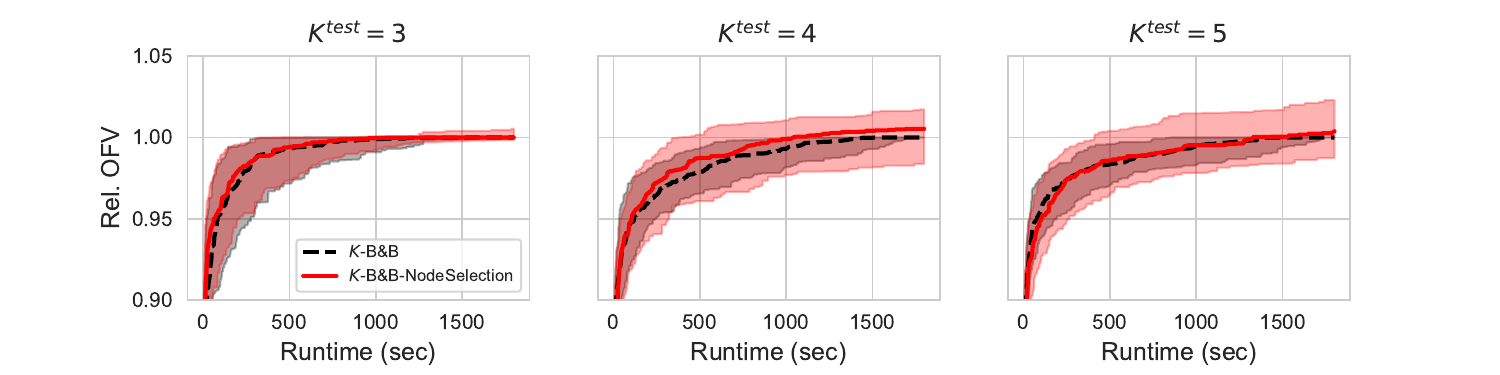} \label{fig:cb_exp6}}
    \caption{EXP5 results of capital budgeting for $K^{test} = K^{train} \in \{3, 4, 5\}$ and EXP6 results for $K^{train} = 6$ and $K^{test} \in \{3, 4, 5\}$ on 100 instances. {The black line gives the average Rel. OFV of \textsc{$K$-B\&B} and the red line that of \textsc{$K$-B\&B-NodeSelection}.}}
\end{figure}

\begin{figure}[htbp!]
    \centering
    \subfloat[EXP5]{
        \includegraphics[width=.9\columnwidth]{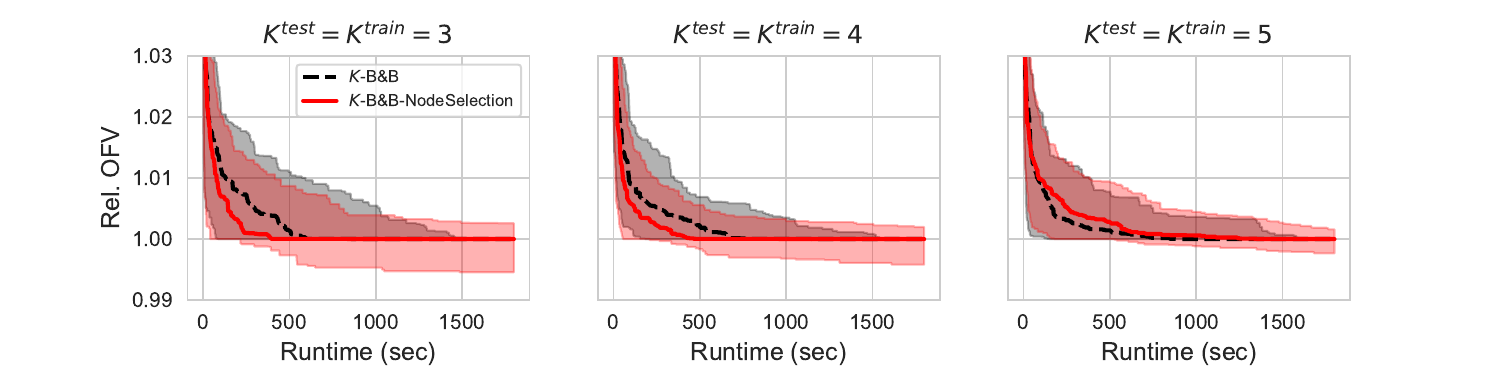} \label{fig:sp_exp5}} \\
    \subfloat[EXP6]{        \includegraphics[width=.9\columnwidth]{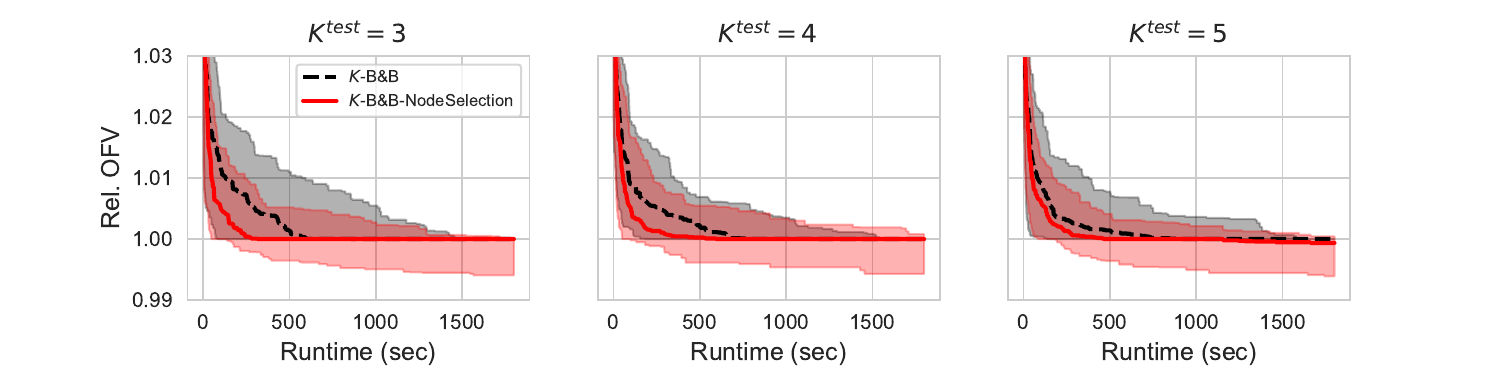} \label{fig:sp_exp6}}
    \caption{EXP5 results of shortest path for $K^{test} = K^{train} \in \{3, 4, 5\}$ and EXP6 results for $K^{train} = 6$ and $K^{test} \in \{3, 4, 5\}$ on 100 instances. {The black line gives the average Rel. OFV of \textsc{$K$-B\&B} and the red line that of \textsc{$K$-B\&B-NodeSelection}.}}
\end{figure}

Now we show the results of the shortest path problem that uses a ML model trained on capital budgeting data. Due to the mismatch of features, we delete the three first-stage-related features not used by shortest path from the capital budgeting data. We can still use the generated capital budgeting data. For an illustration of EXP5, see Figure \ref{fig:sp_exp5}. These plots illustrate that for two out of three values of $K^{test}$, \textsc{$K$-B\&B-NodeSelection} outperforms $K$-B\&B, even though it is trained on data of another problem. More interesting is the following: the performance of different values of $K^{test}$ with $K^{train} = 6$ gives very good results. They are as good as the ones trained on shortest path data. For an illustration of this, see Figure \ref{fig:sp_exp6}.

\subsection{Feature importance} 
\label{sec:experiments_feat_imp}
A trained random forest model allows us to compute feature importance scores. For both problems tested on, these scores are given in Table \ref{tab:feat_imp} for $K^{train} \in \{2, 3, 4, 5, 6\}$. We observe that the relatively highest importance corresponds to the 'Objective' feature (a state feature), with the importance of the other features (both state and scenario ones), being of similar magnitude. Because of that, we cannot draw significant conclusions about the relative importance of the remaining features.
\begin{table}[htbp!]
    \centering    
    \footnotesize
    \caption{Feature importance scores of the selected random forest models for the two problems capital budgeting and shortest path, for each $K^{train}$. The bold scores are given to the scores that belong to the two highest ones of that model. The dashes correspond to the omitted features for the shortest path problem.}
    \label{tab:feat_imp}
\begin{tabular}{ll|ccccc|ccccc}
\toprule
\multirow[c]{2}{*}{\textbf{Feature name}} & & \multicolumn{5}{c|}{\textbf{Capital budgeting}} & \multicolumn{5}{c}{\textbf{Shortest path}} \\ 
 & $K^{train}$ & 2 & 3 & 4 & 5 & 6 & 2 & 3 & 4 & 5 & 6 \\
 \midrule
 Objective                   & &\textbf{ 0.25} & \textbf{0.25} & \textbf{0.24} &\textbf{ 0.28} & \textbf{0.25}   & \textbf{0.10} & \textbf{0.11} & \textbf{0.24} & \textbf{0.13} & \textbf{0.15} \\
Objective difference            & & 0.05 & 0.04 & 0.04 & 0.03 & 0.04                                                & 0.06 & 0.05 & 0.04 & 0.05 & 0.05 \\
Violation                       & & \textbf{0.12} & \textbf{0.10} & 0.07 & 0.07 & 0.06                              & 0.09 & 0.09 & 0.10 & 0.10 & 0.10 \\
Violation difference            & & 0.05 & 0.06 & 0.05 & 0.05 & 0.05                                                & 0.09 & 0.09 & 0.08 & 0.10 & 0.09 \\
Depth                           & & 0.05 & 0.08 & 0.07 & 0.04 & 0.03                                                & \textbf{0.12} & \textbf{0.13} & \textbf{0.12} & \textbf{0.11} & \textbf{0.10} \\
Scenario values                 & & 0.06 & 0.06 & 0.05 & 0.05 & 0.06                                                & 0.09 & \textbf{0.11} & 0.08 & 0.10 & 0.09 \\
Constraint distance             & & 0.08 & 0.06 & 0.05 & 0.05 & 0.05                                                & 0.10 & 0.10 & 0.08 & 0.10 & \textbf{0.10} \\
Scenario distance               & & 0.06 & 0.06 & 0.07 & 0.07 & 0.05                                                & 0.09 & 0.10 & 0.07 & 0.09 & 0.09 \\
Constraint slacks               & & 0.07 & 0.06 & 0.05 & 0.04 & 0.04                                                & 0.09 & 0.06 & 0.05 & 0.07 & 0.07 \\
Det. Objective value            & & 0.05 & 0.06 & \textbf{0.11} & \textbf{0.12} & \textbf{0.12}                     & 0.09 & 0.09 & 0.08 & 0.09 & 0.09 \\
Det. First-stage decisions      & & 0.05 & 0.06 & 0.05 & 0.06 & 0.07                                                & - & - & - & - & - \\
Det. Second-stage decisions     & & 0.00 & 0.00 & 0.00 & 0.00 & 0.00                                                &  0.08 & 0.07 & 0.06 & 0.06 & 0.06 \\
Stat. Objective value           & & 0.06 & 0.06 & 0.08 & 0.10 & \textbf{0.12}                                       & - & - & - & - & - \\
 Stat. Second-stage decisions    & & 0.06 & 0.05 & 0.05 & 0.05 & 0.05                                                & - & - & - & - & - \\
\bottomrule
\end{tabular}
\end{table}

\section{Conclusion and future work} \label{sec:conclusion}
We introduce an ML-based method to improve the $K$-B\&B algorithm \citep{subramanyam2020k} for solving 2SRO optimization problems. $K$-B\&B uses a search tree to optimally partition the uncertainty set into $K$ parts and we propose to use a supervised ML model that learns the best node selection strategy to explore such trees faster. 

For this, we designed a procedure for generating training data and formulated the ML features based on our knowledge of 2SRO so that they are independent of the size, the value of $K$, and the type of problems on which the ML tool is trained. We experimentally show that our method outperforms $K$-B\&B on the problems we test on. We see that when a problem is trained on a smaller instance size, and then applied to the same problem type with bigger instances, our method still outperforms $K$-B\&B, although being less stable. Training and testing on entirely different problem types resulted in mixed results.

As $K$-B\&B has a tree search structure, we believe that our work can be used to tackle other problems solved by a similar tree search structure, wherever expert knowledge can be used to construct meaningful problem size-independent features.

\newpage

\bibliographystyle{plainnat}
\bibliography{bibl.bib}

\newpage
\appendix

\section{Attribute descriptions} \label{app:attributes}
Each scenario gets its own set of attributes. In total there are nine types: one is the scenario vector $\z$ itself ($\a_1^{\z} = \z$), three are determined by the solutions of the main problem, three are extracted from solving the \textit{deterministic problem}, and two are taken from solving the \textit{static problem}.

\paragraph{\textit{Main problem-based attributes}.} For these attributes, only little additional computation is needed. This is due to the fact that the values we use can be taken from the current solution of the main problem. Attributes 2-4 in Table 2 of the paper share this property:

\begin{itemize}
    \item[2.] \textbf{Constraint distance:  } For the first attribute of the `main-problem' type, we look at the constraints generated by the new scenario $\z^*$. We call these constraints the `$\z^*$-constraints'. When a constraint is added to a problem, the resulting feasible region will always be as large as, or smaller than, the feasible region we had before. When the feasible region is large, the objective value we find is often better than for smaller feasible regions (the optimal value found in the large region may be cut off in the small region). Each subset $\barZ_1, \ldots, \barZ_K$ consists of its own set of scenarios, which translates to its own set of constraints in the main problem. These are the `$\barZ_k$-constraints', for all $k \in \K$. These $\barZ_k$-constraints form a feasibility region for the decision pair ($\x, \y_k$). Ideally, we would calculate the volume of the feasible region whenever the $\z^*$-constraints are added to the existing feasible regions. However, obtaining this result is computationally intractable. Therefore we only look at the distance between the constraints. We generate the cosine similarity between the $\z^*$-constraints and the $\barZ_k$-constraints, for each subset separately. When the cosine similarity is high, the distance between the constraints is low. Then, the attribute of the $k$-th child node $\a^{k, \z^*}_2$ is the cosine similarity of the $\z^*$-constraints and the $\barZ_k$-constraints:
    \begin{equation*}
        a^{k, \z^*}_{2, c} = \max_{\z \in \barZ_k} \quad \frac{\gamma_c(\z^*) \cdot \gamma_c(\z)}{\lVert\gamma_c(\z^*)\rVert \lVert \gamma_c(\z)\rVert}, \qquad \forall c \in \{1, \ldots, C\}, \quad \forall k \in \K,
    \end{equation*}
    where $a^{k, \z^*}_{2, c} \in [-1, 1]$ and $\gamma_c: \Z \rightarrow \X \times \Y$ is a function of the left-hand-side of constraint $c \in \{1, \ldots, C\}$ with input scenario $\z$. The first- and second-stage decisions are variable. Finally, the attribute of the $k$-th child node is formulated as $\a^{k, \z^*}_2 = [a^{k, \z^*}_{2, 1}, \ldots, a^{k, \z^*}_{2, C}]$. Thus, the length of  $\a^{k, \z^*}_2$ is equal to the number of uncertain constraints of the MILP formulation of the problem.
    
    \item[3.] \textbf{Scenario distance: } Per subset, we wish to know how far away the new scenario $\z^*$ is from \emph{not} being a violating scenario. We suspect that if the distance is small rather than large, the current solution will not change much. Thus, will not become much worse. This attribute first takes the current solutions of the main problem. Then, for the $k$-th subset it determines the distance between each of the following planes (boundaries of constraints of Eq. (4) in the paper) and the new scenario $\z^* \in \Z$: 
\begin{equation}
\begin{alignedat}{2}
\label{eq:att_3_planes}
    & c(\z)^\intercal \x^* + d(\z)^\intercal \y^*_k - \theta^* = 0, \\
    & T_c(\z)\x^* + W_c(\z)\y^*_k - h_c(\z) = 0, \quad & \forall c \in \{2, \ldots, C\}.
\end{alignedat}
\end{equation}
    Hence, we need to determine the distance between a plane and a point, for each constraint $c \in \{1, \ldots, C\}$, where $c = 1$ corresponds to the objective function.
    The point-to-plane distance of the $c$-th plane for subset $k$ is calculated by projecting $\z^*$ on the normal vector of the plane as follows:
    \begin{equation*}
        \chi^k_c = \frac{|\boldsymbol{\rho}_c(\x^*, \y^k)^\intercal \z^*|}{\lVert\boldsymbol{\rho}_c(\x^*, \y^k)\rVert},
    \end{equation*}
    where $\boldsymbol{\rho}_c: \X \times \Y \rightarrow \Z$ is a vector of coefficients of constraint $c$ of Eq. \eqref{eq:att_3_planes}. To compare the point-to-plane distance of the $K$ subsets, we scale over the sum of the distance of the subsets. Then, the attribute is given as: 
    \begin{equation*}
        a^{k, \z^*}_{3, c} =  \frac{\chi^k_c}{\sum_{k^{'} \in \K} \chi^{k^{'}}_c}, \quad \forall c \in \{1, \ldots, C\}, \quad \forall k \in \K.
    \end{equation*}
    
    \item[4.] \textbf{Constraint slacks: } This attribute takes the slack values of the uncertain constraints of the main problem for the new scenario $\z^*$ with fixed first- and second-stage decisions $\x^*$ and $\y^*$. For the $k$-th subset and the constraints $c \in \{1, \ldots, C\}$, with $c = 1$ the objective constraint, we get:
    \begin{align*}
        s_1^k & = |c(\z^*)^\intercal \x^* + d(\z^*)^\intercal \y^*_k - \theta^*|, \\
        s_c^k & = |T_c(\z^*)\x^* + W_c(\z^*)\y^*_k - h_c(\z^*)|, & \forall c \in \{2, \ldots, C\}.
    \end{align*}
    Similarly as with the previous attribute, we compare the slack values of the subsets by scaling over the sum of the slacks of all subsets: 
    \begin{equation*}
        a^{k, \z^*}_{4, c} =  \frac{s^k_c}{\sum_{k^{'} \in \K} s^{k^{'}}_c}, \quad \forall c \in \{1, \ldots, C\}, \quad \forall k \in \K.
    \end{equation*}
\end{itemize}

\paragraph{\textit{Deterministic problem-based attributes}.}
This problem solves the 2SRO problem for where the newly found scenario $\z^*$ is the only scenario considered. We call this a deterministic problem, since we no longer deal with uncertainty. The problem is formulated as follows:

\begin{equation}
\begin{alignedat}{2}
\min_{\theta^n, \x^n, \y^n} \quad & \theta^n \label{eq:deterministic_problem} \\
    \text{s.t.} \quad & \theta^n \in \R, \x^n\in \X, \y^n \in \Y, \\
    &  \bc(\z^*)^\intercal \x^n + \d(\z^*)^\intercal \y^n \leq \theta^n, \\
    & \T(\z^*)\x^n + \W(\z^*)\y^n \leq \h(\z^*).
\end{alignedat}
\end{equation}

\noindent By solving this problem we obtain Attributes 5-7: 
\begin{itemize}
    \item[5.] \textbf{Deterministic objective} function value $\theta^n$,
    \item[6.] \textbf{Deterministic first-stage} decisions $\x^n$,
    \item[7.] \textbf{Deterministic second-stage} decisions $\y^n$.
\end{itemize}

\paragraph{\textit{Static problem-based attributes}.} 

A very simple method for approximately solving the 2SRO problem is to first solve the first-stage decision for all scenarios in the uncertainty set. Then, after a realization of uncertainty, we combine this first-stage decision with the scenario to determine the second-stage decision. This is an naive way of solving 2SRO, thus not optimal for all scenarios. But, it could give us some information on the approximate solutions to scenarios in the problem. Solving the static problem consists of two steps: First, we obtain the static robust first-stage decisions $\Bar{\x}$ by solving

\begin{equation}
\begin{alignedat}{2}
\min_{\theta, \x, \y} \quad & \theta \label{eq:static_prep} \\
    \text{s.t.} \quad & \theta \in \R, \x \in \X, \y \in \Y,\\
    &  \bc(\z)^\intercal \x + \d(\z)^\intercal \y \leq \theta, \qquad \qquad & \forall \z \in \Z, \\
    & \T(\z)\x + \W(\z)\y \leq \h(\z), & \forall \z \in \Z.
\end{alignedat}
\end{equation}
This problem can be reformulated via the mathematically tractable formulation presented in \cite{ben2009robust}. Secondly, by fixing $\x$ to $\Bar{\x}$, we obtain the objective value $\theta^s$ and second-stage decisions $\y^s$ by solving
\begin{equation}
\begin{alignedat}{2}
\min_{\theta^s, \y^s} \quad & \theta^s \label{eq:static_problem} \\
    \text{s.t.} \quad & \theta^s \in \R, \y^s \in \Y,\\
    &  \bc(\z^*)^\intercal \Bar{\x} + \d(\z^*)^\intercal \y^s \leq \theta^s, \\
    & \T(\z^*) \Bar{\x} + \W(\z^*)\y^s \leq \h(\z^*).
\end{alignedat}
\end{equation}
Note that problem \eqref{eq:static_prep} is solved only once in the algorithm, while problem \eqref{eq:static_problem} needs to be solved for each scenario. By solving this problem we obtain Attributes 8 and 9: 
\begin{enumerate}
    \item[8.] \textbf{Static objective} function value $\theta^s$,
    \item[9.] \textbf{Static second-stage}  decisions $\y^s$ 
\end{enumerate}

\clearpage
\section{Omitted pseudocodes} \label{app:pseudocodes}
The $K$-adaptability branch-and-bound ($K$-B\&B) algorithm, as presented in \cite{subramanyam2020k}, is given in Algorithm \ref{algorithm:k_bb}. In our implementation of this algorithm, we apply a depth-first search strategy and select a random node of $\N$ instead of the first one.

\begin{algorithm}[H]
\footnotesize
\caption{\textsc{$K$-B\&B}}
\label{algorithm:k_bb}
\DontPrintSemicolon
\KwInput{\quad Problem instance $\P(N)$ with size $N$, \\ \quad number of partitions $K$}
\KwOutput{\quad Objective value $\theta$, first-stage decisions $\x$, second-stage decisions $\y = \{\y_1, \ldots, \y_K \}$, \\
\quad subsets with scenarios $\barZ_k$ for all $k \in \{1,\ldots,K\}$}
\KwInit{\quad Incumbent partition: $\tau^i := \{\barZ_1, \ldots, \barZ_K\}$, where $\barZ_k = \emptyset$ for all $k \in \K$, \\
        \quad set containing all node partitions yet to explore: $\N := \{\tau^i\}$, \\
        \quad incumbent solutions: $(\theta^i, \x^i, \y^i) := (\infty, \emptyset, \emptyset)$
}
\BlankLine
\While{$\N$ not empty}{
    Select the first node with partition $\tau= \{\barZ_1, \ldots, \barZ_K \}$ from $\mathcal{N}$, then $\mathcal{N} \leftarrow \mathcal{N}\setminus \{\tau\}$ \label{alg_step:select_node}\\ 
    
    $(\theta^*, \x^*, \y^*) \leftarrow \textit{main problem}(\tau)$ \\
    \If{$\theta^* > \theta^i$}{
        Prune tree since current objective is worse than best solution found. \label{alg_step:prune}\\ 
        Continue to line \ref{alg_step:select_node}.}
    $(\z^*, \zeta^*) \leftarrow \textit{subproblem}(\theta^*, \x^*, \y^*)$ \\
    \uIf{$\zeta^* > 0$}{
        Solution not robust, create $K$ new branches. \\
        \For{$k \in \{1, \ldots, K\}$}{
            $\tau^k := \{\barZ_1, \ldots, \barZ_k \cup \{\z^*\}, \ldots, \barZ_K\}$ \\
            $\N \leftarrow \N \cup \{\tau^k\} $ \label{alg_step:branch}
            }
        }
    \Else{
        Current solution robust, prune tree. \label{alg_step:robust}\\
        $(\theta^i, \x^i, \y^i, \tau^i) \leftarrow (\theta^*, \x^*, \y^*, \tau)$
        }
    }
    \Return{$(\theta^i, \x^i, \y^i, \tau^i)$}
\end{algorithm}

The steps for obtaining the ML model used for node selection consists of two parts: (i) making training data and (ii) training the ML model. More details on these steps are given in Procedure \ref{proc:strategy_model}.

\renewcommand*{\algorithmcfname}{Procedure}
\begin{algorithm}[H]
\footnotesize
\caption{\textsc{StrategyModel}}
\label{proc:strategy_model}
\DontPrintSemicolon
\KwInput{\quad Train instances $\P_1^{train}(N^{train}), \ldots, \P_I^{train}(N^{train})$ \\ \\ \quad number of partitions for training $K^{train}$, \\ \quad level for training $L^{train}$, \\ \quad quality threshold $\epsilon$, \\
\quad $R$ for random dives per node}
\KwOutput{\quad Trained node selection strategy \textit{model}.}
\BlankLine
\tcp{Get training data for $I$ instances}
\For{$i \in \{1, \ldots, I\}$}
    {$(D, \p)_i \leftarrow \textit{generate train data}(\P^{train}_i(N^{train}), K^{train}, L^{train}, R)$
    $\q_i \leftarrow \textit{quality}(\p_i, \epsilon)$
    }
\BlankLine
\tcp{Train node selection strategy model}
Set $\textit{model} \leftarrow \textit{train ML model}(\{(D, \q)_1, \ldots, (D, \q)_{I}\})$ \\
\BlankLine
\Return{\textit{model}}
\end{algorithm}
For scaling some of the features (see Section 3.2 of the paper), information needs to be gathered by performing several initial dives in the tree. The steps of these dives are given in Procedure \ref{proc:init_run}.

\clearpage

\begin{algorithm}[H]
\footnotesize
\caption{\textsc{InitialRun}}
\label{proc:init_run}
\DontPrintSemicolon
\KwInput{\quad Test instance $\P^{test}(N^{test})$, \\ \quad number of partitions for testing $K^{test}$}
\KwOutput{\quad Objective value $\theta$, first-stage decisions $\x$, second-stage decisions $\y = \{\y_1, \ldots, \y_K \}$, \\
\quad Subsets with scenarios $\barZ_k$ for all $k \in \{1,\ldots,K\}$}
\KwInit{\quad Incumbent partition: $\tau^i := \{\barZ_1, \ldots, \barZ_K\}$, where $\barZ_k = \emptyset$ for all $k \in \K$, \\}
\BlankLine
Set $\tau = \tau^i$ \\
\While{solution not robust}{
    $(\theta^*, \x^*, \y^*) \leftarrow \textit{main problem}(\tau)$ \\
    $(\z^*, \zeta^*) \leftarrow \textit{subproblem}(\theta^*, \x^*, \y^*)$ \\
    \uIf{$\zeta^* > 0$}{
        Solution not robust, random node selection. \\
        $k^{'} \leftarrow \textit{random uniform sample}([1, K])$ \\
        $\tau := \{\barZ_1, \ldots, \barZ_{k^{'}} \cup \{\z^*\}, \ldots, \barZ_K\}$ \\
        }
    \Else{
        Robust solution found \\
        $\textit{scaling info} \leftarrow (\theta^0, \zeta^0, \kappa^0)$
        }
    }
\Return{\textit{scaling info}}
\end{algorithm}

\section{Problem formulations} \label{app:problems}
In the experiments, our method has been tested for several problems. The descriptions of these problems and their MILP formulations are given in this section.
\subsection{Capital budgeting with loans} \label{app:problems_cb}
We consider the capital budgeting with loans problem as defined in \cite{subramanyam2020k}, where a company wishes to invest in a subset of $N$ projects. Each project $i$ has an uncertain cost $c_i(\z)$ and an uncertain profit $r_i(\z)$, defined as
\begin{equation*}
    c_i(\z) = \big (1 + \boldsymbol{\Phi}_i^\intercal \z/2 \big ) c_i^0 \quad \text{and} \quad r_i(\z) = \big (1 + \boldsymbol{\Psi}_i^\intercal \z/2 \big ) r_i^0, \quad \forall i \in \{1,\ldots,N\},
\end{equation*}
where $c_i^0$ and $r_i^0$ represent the nominal cost and the nominal profit of project $i$, respectively. $\boldsymbol{\Phi}_i^\intercal$ and $\boldsymbol{\Psi}_i^\intercal$ represent the $i$-th row vectors of the sensitivity matrices $\boldsymbol{\Phi}, \boldsymbol{\Psi} \in \R^{N \times N_z}$.
The realizations of the uncertain vector $\z$ belong to the uncertainty set $\Z = [-1, 1]^{N_z}$, where $N_z$ is the dimension of the uncertainty set.  

The company can invest in a project either before or after observing the risk factor $\z$. In the latter case, the company generates only a fraction $\eta$ of the profit, which reflects a penalty of postponement. However, the cost remains the same as in the case of an early investment. The company has a given budget $B$, which the company can increase by loaning from the bank at a unit cost of $\lambda >0$, before the risk factors $\z$ are observed. A loan after the observation occurs, has a unit cost of $\mu \lambda$, with $\mu > 1$. The objective of the capital budgeting problem is to maximize the total revenue subject to the budget. This problem can be formulated as an instance of the $K$-adaptability problem as follows:
\begin{align*}
    \max_{(x_0,\x) \in \X, (y_0,\y) \in \Y^K} \min_{\z \in \Z} \max_{k \in \K} \quad & \theta \\ 
    \text{s.t.} \quad  & r(\z)^\intercal(\x + \eta\y_k) - \lambda (x_0 + \mu y_0^k) \geq \theta, \\
    & \x + \y_k \leq \boldsymbol{e}, \\ 
    & \boldsymbol{c(z)}^\intercal \x \leq B + x_0, \\
    & \boldsymbol{c(z)}^\intercal (\x + \y_k) \leq B + x_0 + y_0^k, 
\end{align*}
where $\X = \Y = \R_+ \times \{0, 1\}^N$, $y_0 = \{y_0^1, \ldots, y_0^K\}$, $\y = \{\y_1, \ldots, \y_K\}$, $x_0$ and $y_0$ are the amounts of taken loan in the first and second stage, respectively. Moreover, $x_i$ and $y_i$ are the binary variables that indicate whether we invest in the $i$-th project in the first- and second-stage, respectively. The constraints $\boldsymbol{c(z)}^\intercal \x \leq B + x_0$ ensure that for the first stage, the expenditures are not more than the budget plus the loan taken before the realization of uncertainty.

\paragraph{Test case.}
Similarly as in \cite{subramanyam2020k}, the uncertainty set dimension $N_z$ is set to $N_z = 4$. The nominal cost vector $\boldsymbol{c}^0$ is chosen uniformly at random from the set $[0, 10]^N$. Let $\boldsymbol{r}^0 = \boldsymbol{c}^0/5$, $B = \boldsymbol{e}^\intercal \boldsymbol{c}^0/2$, and $\eta=0.8$. The rows of the sensitivity matrices $\boldsymbol{\Phi}$ and $\boldsymbol{\Psi}$ are sampled uniformly from the $i$-th row vector, which is sampled from $[0, 1]^{N_z}$, such that $\boldsymbol{\Phi}_i^\intercal \boldsymbol{e} = \boldsymbol{\Psi}_i^\intercal \boldsymbol{e} = 1$ for all $i \in \{1,\ldots,N\}$. This is also known as the unit simplex in $\R^{N_p}$. For determining the cost of the loans, we set $\lambda = 0.12$ and $\mu = 1.2$.

\subsection{Shortest path} \label{app:problems_sp}
We consider the shortest path problem with uncertain arc weights as defined in \cite{subramanyam2020k}. Let $G=(V,A)$ be a directed graph with nodes $V = \{1, ..., N\}$, arcs $A \subseteq V \times V$, and arc weights $d_{ij}(\z) = (1 + z_{ij}/2)d_{ij}^0$, $(i,j) \in A$. Where $d_{ij}^0 \in \R_+$ represents the nominal weight of the arc $(i, j) \in A$ and $z_{ij}$ denotes the uncertain deviation from the nominal weight. The uncertainty set is defined as
\begin{equation*}
    \Z = \Big \{\z \in [0, 1]^{|A|}: \sum_{(i,j) \in A} z_{ij} \leq \Gamma \Big \}.
\end{equation*}
This uncertainty set imposes that at most $\Gamma$ arc weights may maximally deviate from their nominal values. We need to find the shortest path from the source node $s$, to the sink node $t$ before observing the realized arc weights. This shortest problem can be formulated as an instance of the $K$-adaptability problem
\begin{align*}
    \min_{\y \in \Y^K} \max_{\z \in \Z} \min_{k \in \K} \quad & \theta \\
    \text{s.t.} \quad & {\d(\z)}^\intercal \y_k \leq \theta,  \\
     & \sum_{(j,l) \in A} y^k_{jl} - \sum_{(i,j) \in A} y^k_{ij} \geq \1_{\{j=s\}} - \1_{\{j=t\}}, \qquad \forall j \in V, \\
    & \Y \subseteq \{0, 1\}^{|A|}.
\end{align*}
Note that this problem contains only binary second-stage decisions and uncertainty in the objective function. The $K$-B\&B algorithm, will find $K$ shortest paths from $s$ to $t$. After $\z$ is observed, the path $\y_k$ will be chosen if $\z \in \Z_k$.

\paragraph{Normal test case.}
The coordinates in $\mathbb{R}^2$ for each vertex $i \in V$ are uniformly chosen at random from the region $[0, 10]^2$. The nominal weight of the arc $(i, j) \in A$ is the Euclidean distance between node $i$ and $j$. The source node $s$ and the sink node $t$ are defined to be the nodes with the maximum nominal distance between them. The $\lfloor{0.9(N^2-N)}\rfloor$ arcs with the highest nominal weight will be deleted to define the arc set $A$. The budget of the uncertainty set $\Gamma$ is set to seven.

\paragraph{Sphere test case.}
The instances of this type have nodes that are spread over a sphere. This is done as follows. First, each node in the three-dimensional graph is sampled from the standard normal distribution and then normalized. The distance between node $i$ and $j$ is then derived by its spherical distance. To obtain this, first the Euclidean distance $d_{ij}$ between node $i$ and $j$ is computed. Then, the arc sine of $d_{ij}/2$ is computed to get the spherical distance. The $\lfloor{0.7(N^2-N)}\rfloor$ arcs with the highest nominal weight will be deleted to define the arc set $A$. The budget of the uncertainty $\Gamma$ is set to seven.

\subsection{Knapsack}
We consider the two-stage version of the knapsack problem where the profit per item is uncertain. This formulation is based on that of \cite{buchheim2017min}. Let $N$ be the number of items, $p_i(\z) = (1-z_i/2)p^0_i$ the profit for item $i \in \{1, \ldots, N\}$, where $p^0_i \in \R$ is the nominal profit value, $z_i$ is the deviation, $\w \in \R^N$ is the weight vector, and $C = c\sum_{i=1}^N w_i$ is the total capacity of the knapsack with $c \in (0, 1)$. The uncertainty set is defined as
\begin{equation*}
    \Z = \Big \{\z \in [0, 1]^{N}: \sum_{i=1}^N z_i \leq \Gamma \Big \},
\end{equation*}
where $\Gamma = \gamma N$ and $\gamma \in (0, 1)$. This problem can be formulated as an instance of the $K$-adaptability problem
\begin{align*}
    \max_{\y \in \Y^K} \min_{\z \in \Z} \max_{k \in \K} \quad & \theta \\
    \text{s.t.} \quad & {\p(\z)}^\intercal \y_k \geq \theta, \\
    & \w^\intercal \y_k \leq C, \\
    & \Y \subseteq \{0, 1\}^{N},
\end{align*}
where $y_i$ is the decision of putting item $i$ in the knapsack.

\paragraph{Test case.}
The weight $w_i$ of each item $i \in \{1, \ldots, N\}$ is uniformly chosen at random from $[1, 15]$ and the cost $c_i$ from $[100, 150]$. The values of $c$ and $\gamma$ are selected in the experiments section of the paper.

\section{Parameter tuning}
For both training the ML model (a random forest) and applying it to a problem, we have defined multiple parameters in \textsc{StrategyModel} and \textsc{$K$-B\&B-NodeSelection}. The tuning of these parameters is explained in this section.

\subsection{Capital budgeting with loans} \label{app:param_tuning_cb}
For each $K \in \{2, \ldots, 6\}$ we have trained five random forests: each for $\epsilon \in \{0.05, 0.1, 0.2, 0.3, 0.4\}$. We have first decided on the values of $\iota$ for generating training data. Problems become more complex when $K$ grows, which results in more time needed per dive. For tuning these parameters, we have generated training data by running the algorithm for two hours. Thus, we have fixed the parameter $T$ to two. In Table \ref{tab:cb_data}, for each combination of $K$ and $\iota \in \{2, 5, 10, 15\}$ the following information is shown: number of data points, number of searched instances $I$, and the average of $L^{train}$ reached per instance. 

\begin{table}[htbp!]
    \centering
    \caption{Generated training data info for combinations of $K$ and $\iota$. (num. data points, $I$, average $L^{train}$).}
    \label{tab:cb_data}
    \footnotesize
    \begin{tabular}{c|cccc}
    \toprule
    \multirow[c]{2}{*}{$K$} & \multicolumn{4}{c}{$\iota$ (in minutes)} \\ 
    {} &     2  &              5  &              10 &             15 \\
    \midrule
    2 &  (17284, 60, 9) &  (7192, 24, 10) &             - &            -     \\
    3 &  (37227, 60, 7) &  (36510, 24, 8) &  (34536, 12, 9) &  (30510, 8, 9) \\
    4 &  (21572, 60, 5) &  (21860, 24, 6) &  (24312, 12, 7) &  (17008, 8, 7) \\
    5 &  (22425, 60, 5) &  (23830, 24, 5) &  (22510, 12, 6) &  (20040, 8, 6) \\
    6 &  (17820, 60, 4) &  (17544, 24, 5) &  (18834, 12, 5) &  (17382, 8, 5) \\
    \bottomrule
    \end{tabular}
\end{table}
We have then selected per $K$ a value of $\iota$ that has high values of the number of data points, $I$ and $L^{train}$. Then, for each $\epsilon$ and $K$, the dataset related to this value of $\iota$ has been trained. The accuracy of these models are given in Table \ref{tab:cb_accuracy}. 

\begin{table}[htbp!]
    \centering
    \footnotesize
    \caption{Test accuracy of random forest for combinations of $K$ (with best $\iota$) and the threshold $\epsilon$.}
    \label{tab:cb_accuracy}
    \begin{tabular}{l|lllll}
    \toprule
    \multirow[c]{2}{*}{$K$ $(\iota)$} & \multicolumn{5}{c}{$\epsilon$} \\ 
     {}    &    0.05  &    0.1 &    0.2 &    0.3 &    0.4 \\
    \midrule
    2 (2)  &    0.971 &  0.988 &  0.971 &  0.988 &  0.983 \\
    3 (5)  &    0.929 &  0.959 &  0.959 &  0.967 &  0.981 \\
    4 (5)  &    0.922 &  0.950 &  0.977 &  0.982 &  0.995 \\
    5 (5)  &    0.958 &  0.937 &  0.967 &  0.975 &  0.992 \\
    6 (10) &    0.937 &  0.952 &  0.968 &  0.974 &  0.984 \\
    \bottomrule
    \end{tabular}
\end{table}

\begin{figure}[htbp!]
    \centering
    \includegraphics[width=0.6\columnwidth]{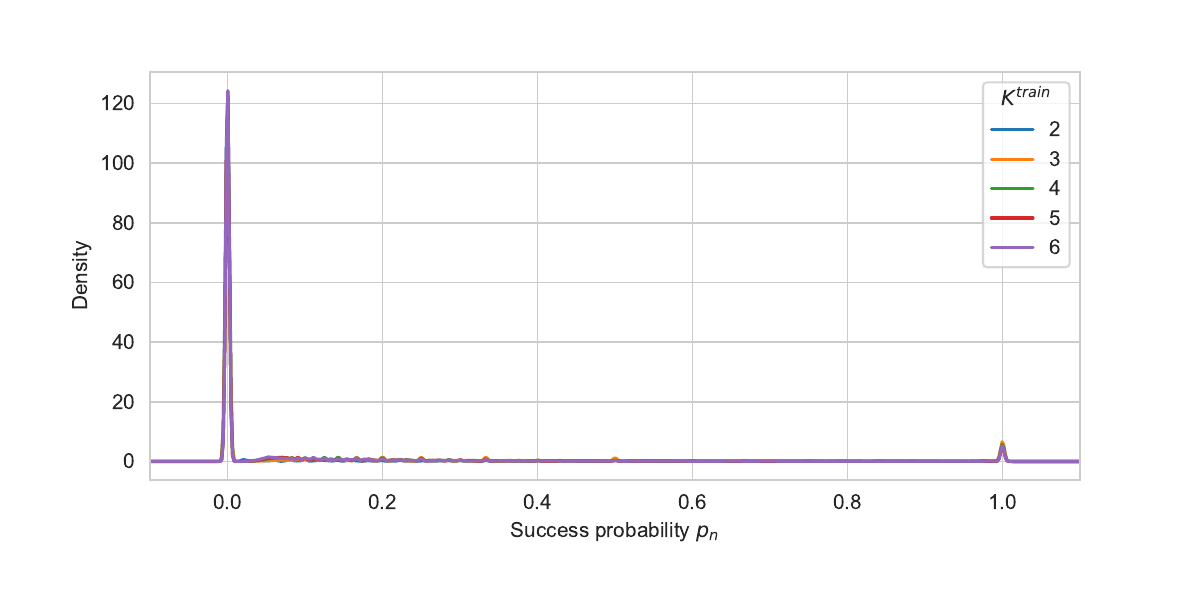}
    \caption{Density of the success probability $p_n$ for each value of $K^{train}$.}
    \label{fig:cb_suc_pred_dens}
\end{figure}

The table above shows that $\epsilon$ does not influence the accuracy of the model. However, Figure 9 shows that the algorithm performs better when $\epsilon$ is very small. If we look at the density of success probabilities $p$ in Figure \ref{fig:cb_suc_pred_dens}, we notice that the vast majority of data points have $p_n \approx 0$. These two observations indicate that any value of $p_n$ slightly higher than zero is special, and the corresponding node is considered as a good node to visit. 

\begin{figure}[htbp!]
    \centering
    \includegraphics[width=\columnwidth]{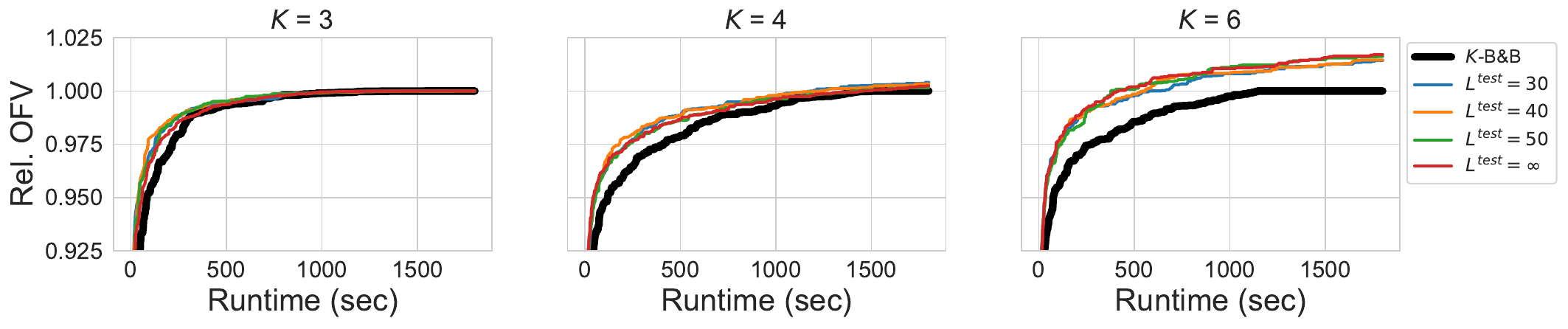}
    \caption{Results of $K$-B\&B with random dives and \textsc{$K$-B\&B-NodeSelection} with combinations of $K$ and bigger values of $L^{test}$.}
    \label{fig:cb_Ldiff_big}
\end{figure}

In the experiments section of the paper, we noticed that high values of $L^{test}$ outperformed lower ones in the \textsc{$K$-B\&B-NodeSelection} algorithm. In Figure \ref{fig:cb_Ldiff_big}, we show the results for higher values of $L^{train}$ than the ones shown in Figure 9 in the paper for a fixed $\epsilon = 0.05$.

\subsection{Shortest path on a sphere} \label{app:param_tuning_sp}
For the shortest path problem we also want to decide on the parameter $\iota$ per $K \in \{2, 3, 4, 5, 6\}$. We noticed that the total number of scenarios needed until a robust solution is found, is larger for shortest path than for capital budgeting. Therefore, the duration per training instance should increase. The range of the number of minutes is $\iota \in \{5, 10, 15, 20\}$. In Table \ref{tab:sp_data}, for each combination of $K$ and $\iota$ the number of data points, instances, and training level $L^{train}$ is given. For now, $T = 2$ is fixed. 

\begin{table}[htbp!]
    \centering
    \caption{Generated training data info for combinations of $K$ and $\iota$. (num. data points, $I$, average $L^{train}$).}
    \label{tab:sp_data}
    \footnotesize
    \begin{tabular}{c|cccc}
    \toprule
    \multirow[c]{2}{*}{$K$} & \multicolumn{4}{c}{$\iota$ (in minutes)} \\ 
    {} &          5  &             10  &             15 &            20  \\
    \midrule
    2 &  (10802, 24, 8) &   (8290, 12, 9) &  (9160, 8, 10) &  (8494, 6, 10) \\
    3 &  (8370, 24, 6) &   (6726, 12, 6) &   (9858, 8, 6) &   (7506, 6, 7) \\
    4 &  (9884, 24, 5) &   (8496, 12, 6) &   (7916, 8, 6) &  (15872, 6, 6) \\
    5 &  (13795, 24, 4) &   (7415, 12, 5) &  (12490, 8, 5) &  (21110, 6, 6) \\
    6 &  (26346, 24, 5) &  (10794, 12, 5) &  (18948, 8, 5) &  (25788, 6, 5) \\
    \bottomrule
    \end{tabular}
\end{table}
We have then selected per $K$ a value of $\iota$ that has high values of the number of data points, $I$ and $L^{train}$. Then, for each $\epsilon$ and $K$, the dataset related to this value of $\iota$ has been trained. See Table \ref{tab:sp_accuracy} for the accuracy of these models.
\begin{table}[htbp!]
    \centering
    \footnotesize
    \caption{Number of data points used for training and test accuracy of random forest for combinations of $K$ (with best $\iota$) and hours spent for generating training data $T$, given the threshold $\epsilon = 0.05$. (num. data points, test accuracy).}
    \label{tab:sp_accuracy}
    \begin{tabular}{l|llll}
    \toprule
        \multirow[c]{2}{*}{$K (\iota)$} & \multicolumn{4}{c}{$T$}    \\ 
        {}     &  1    &    2   &   5    &    10    \\
        \midrule
        2 (15) & (6586, 0.955) &   (9160, 0.946) &  (16852, 0.923) &  (35674, 0.947) \\
        3 (15) & (4170, 0.952) &   (9858, 0.939) &  (28347, 0.919) &  (47346, 0.937) \\
        4 (20) & (2828, 0.931) &  (15872, 0.969) &  (37652, 0.966) &  (69244, 0.932) \\
        5 (20) & (4790, 0.875) &  (21110, 0.972) &   (47000, 0.93) &   (77735, 0.91) \\
        6 (15) & (13500, 0.948) &  (18948, 0.916) &  (54564, 0.945) &  (91254, 0.955) \\
    \bottomrule
    \end{tabular}
\end{table}

We noticed for the shortest path problem that more data points significantly increased the performance of the ML model on the algorithm. An overview of the chosen values of $\iota$ and $T$ (hours spent for getting training data) per $K$ are given in Table \ref{tab:sp_params}.

\begin{table}[htbp!]
    \centering
    \footnotesize
    \caption{Chosen parameter combination for each $K$. The values of $\epsilon$ and $L^{test}$ are fixed to 0.05 and $\infty$, respectively.}
    \label{tab:sp_params}
    \begin{tabular}{lll}
    \toprule
        $K$ & $\iota$ & $T$       \\ 
        \midrule
        2 & 15 & 10  \\
        3 & 15 & 5  \\
        4 & 20 & 5  \\
        5 & 20 & 10  \\
        6 & 15 & 10  \\
    \bottomrule
    \end{tabular}
\end{table}

\begin{figure}[htbp!]
    \centering
     \includegraphics[width=0.6\columnwidth]{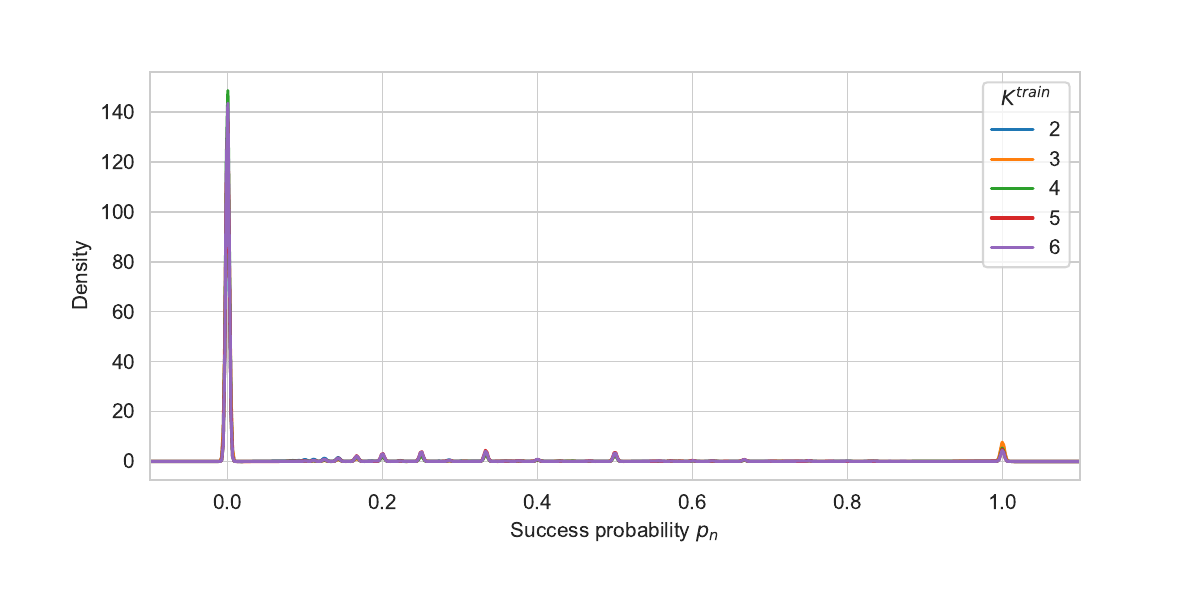}
    \caption{Density of the success probability $p_n$ for each value of $K^{train}$.}
    \label{fig:sp_suc_pred_dens}
\end{figure}

The density of the success probabilities for the data points of the shortest path problem is given in Figure \ref{fig:sp_suc_pred_dens}. This is very similar to the density of capital budgeting (see Figure \ref{fig:cb_suc_pred_dens}).

\clearpage

\section{Full Results} \label{app:results}
We have applied \textsc{$K$-B\&B-NodeSelection} to multiple problems, where the training and testing instance specifications also varied. In the main body, only a subset of the experiments are shown. In this section, all of them are given.

\subsection{Capital budgeting with loans} \label{app:results_cb}
\begin{figure}[htbp!]
\centering
    \includegraphics[width=\columnwidth]{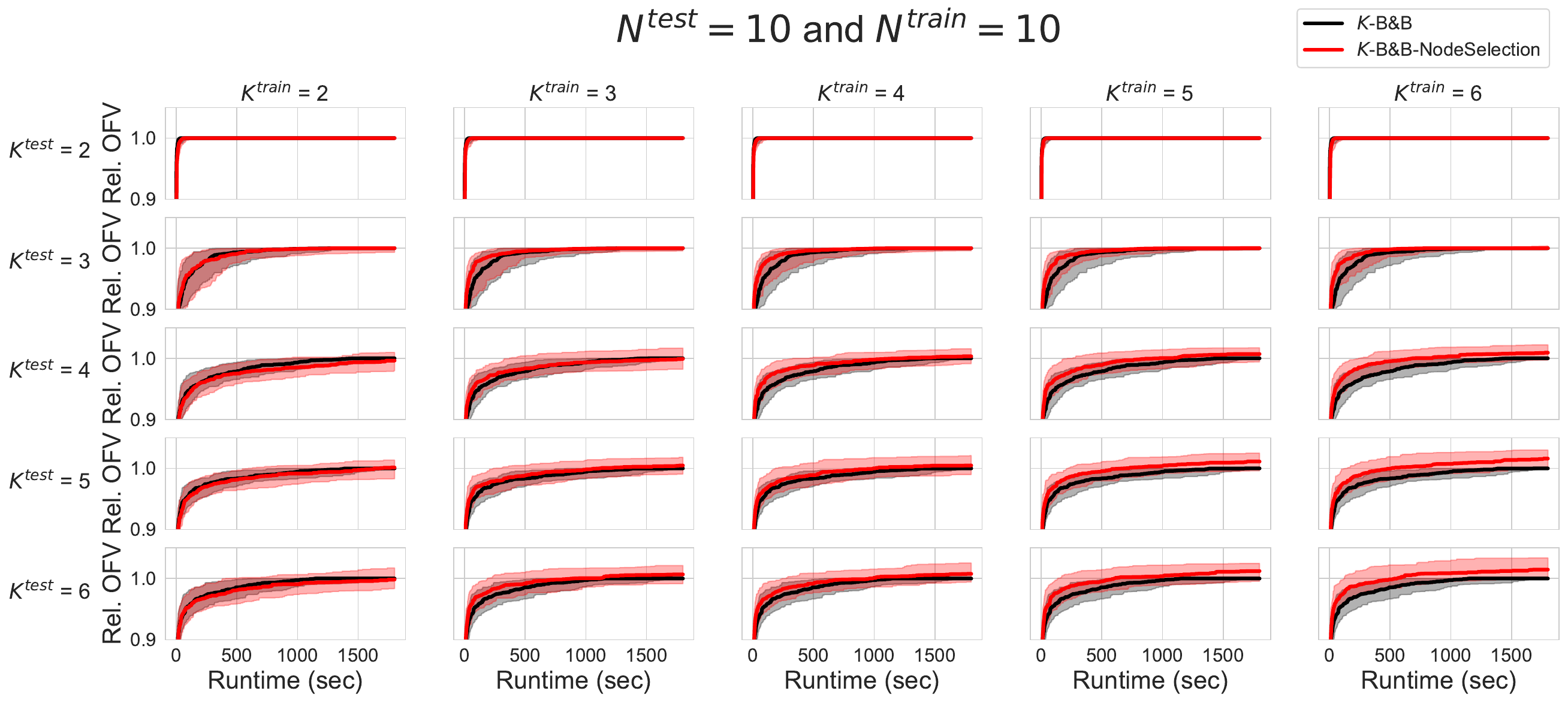}
        \caption{Comparison of results between $K$-B\&B and \textsc{$K$-B\&B-NodeSelection} for 100 instances of the capital budgeting problem. The results of EXP1 and EXP2 are shown, where $N^{test}=N^{train}=10$. The regions with shaded color around the curves denote its $75\%$ CI.}
        \label{fig:cb_train_test_diffK_10}
\end{figure}

\begin{figure}[H]
\centering
    \includegraphics[width=\columnwidth]{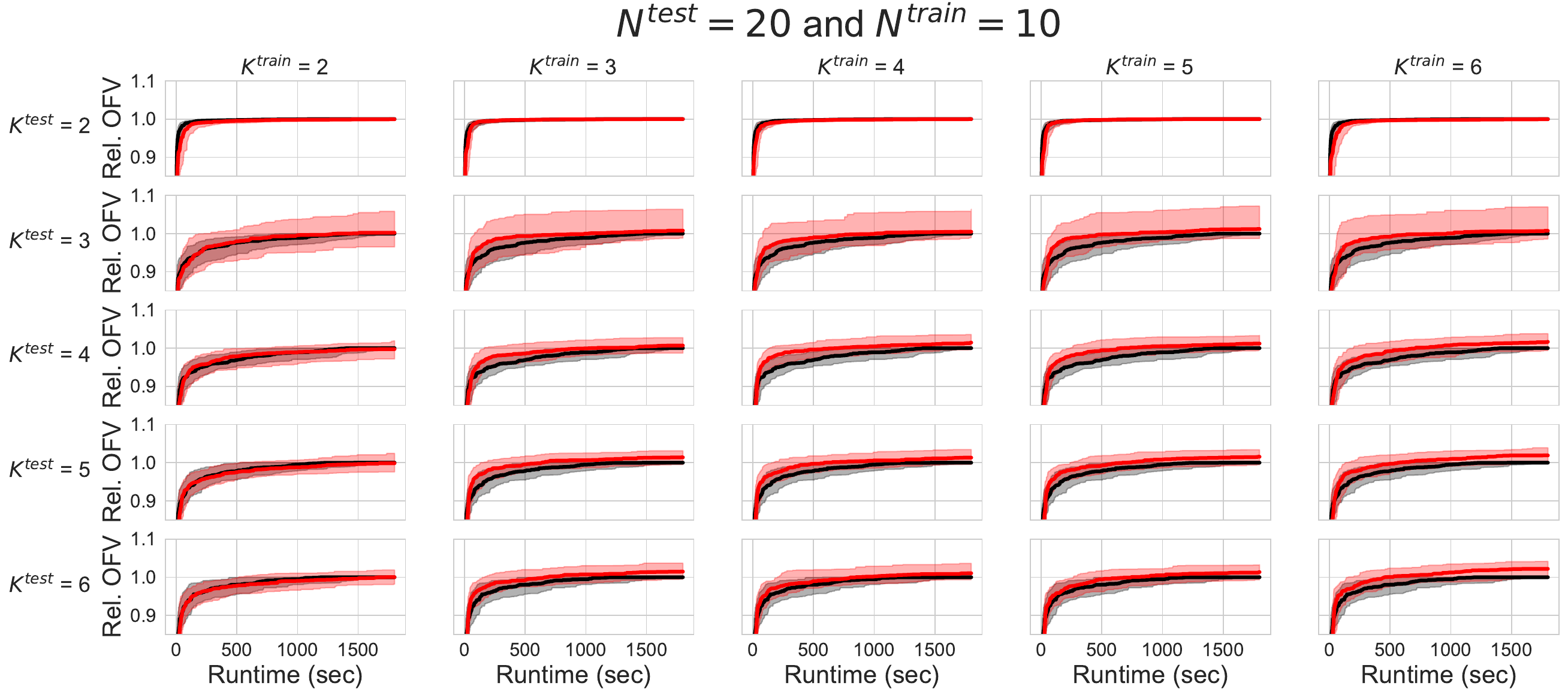}
        \caption{Comparison of results between $K$-B\&B and \textsc{$K$-B\&B-NodeSelection} for 100 instances of the capital budgeting problem. The results of EXP3 and EXP4 are shown, where $N^{train}=10, N^{test}=20$. The regions with shaded color around the curves denote its $75\%$ CI.}
        \label{fig:cb_train_test_diffK_20}
\end{figure}

\begin{figure}[H]
\centering
    \includegraphics[width=\columnwidth]{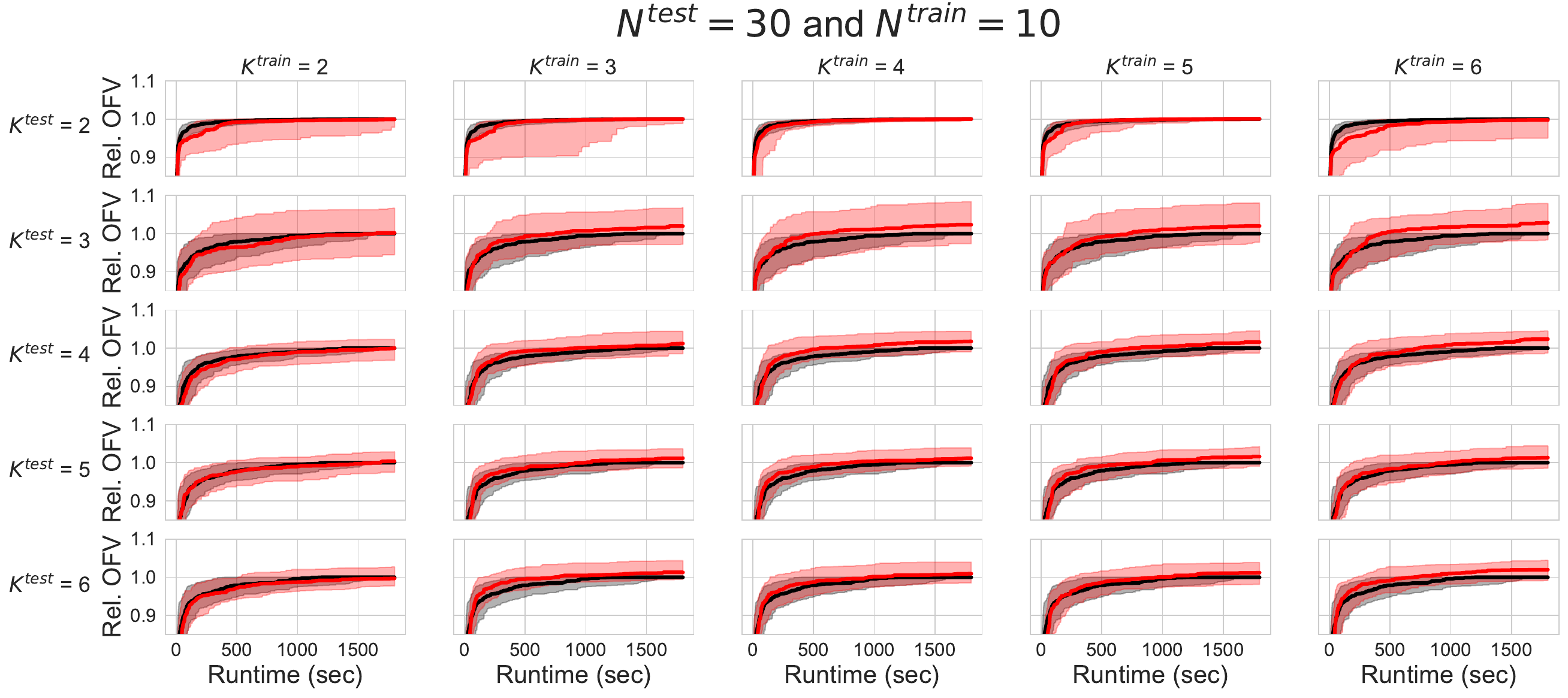}
        \caption{Comparison of results between $K$-B\&B and \textsc{$K$-B\&B-NodeSelection} for 100 instances of the capital budgeting problem. The results of EXP3 and EXP4 are shown, where $N^{train}=10, N^{test}=30$. The regions with shaded color around the curves denote its $75\%$ CI.}
        \label{fig:cb_train_test_diffK_30}
\end{figure}

\subsection{Shortest path} \label{app:results_sp}
\begin{figure}[H]
\centering
    \includegraphics[width=\columnwidth]{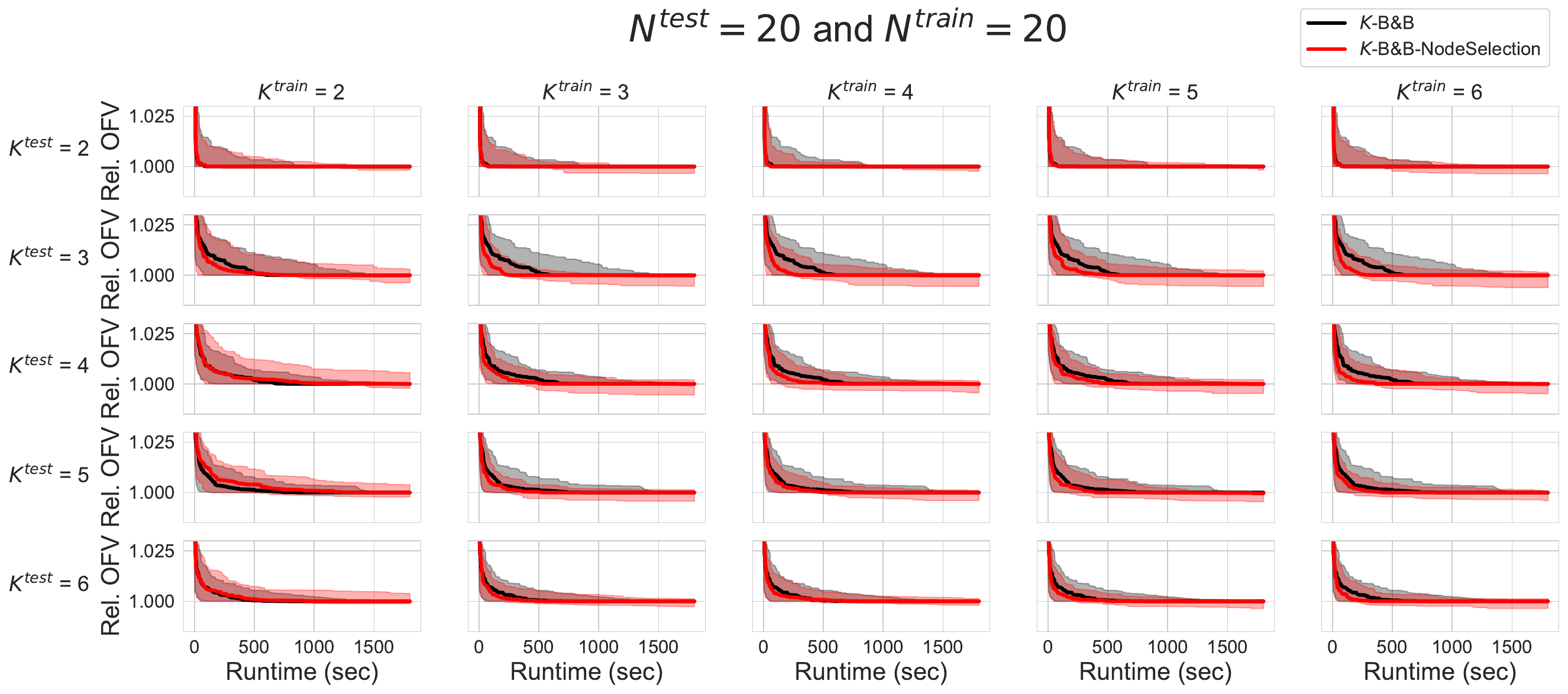}
        \caption{Comparison of results between $K$-B\&B and \textsc{$K$-B\&B-NodeSelection} for 100 instances of the shortest path problem. The results of EXP1 and EXP2 are shown, where $N^{test}=N^{train}=20$. The regions with shaded color around the curves denote its $75\%$ CI.}
        \label{fig:sp_train_test_diffK_20}
\end{figure}

\begin{figure}[H]
\centering
    \includegraphics[width=\columnwidth]{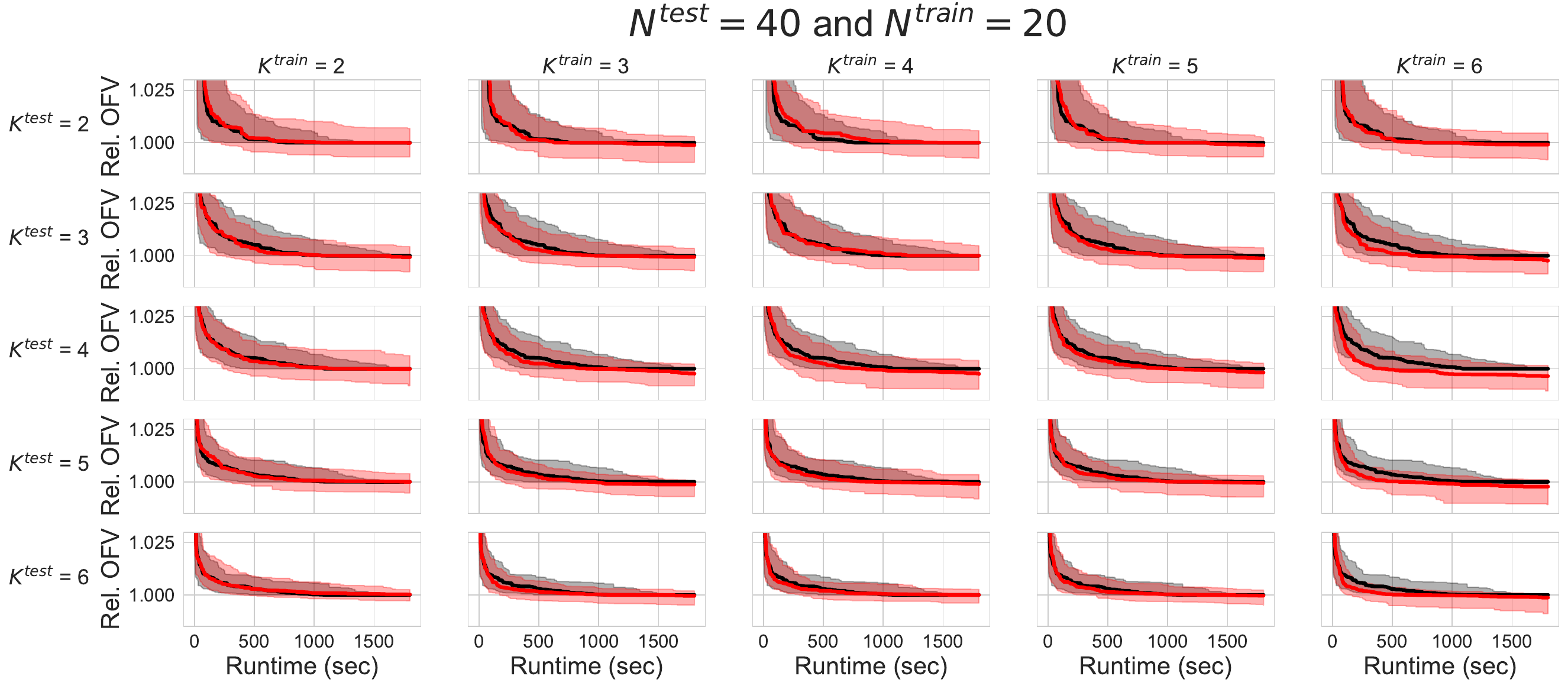}
        \caption{Comparison of results between $K$-B\&B and \textsc{$K$-B\&B-NodeSelection} for 100 instances of the shortest path problem. The results of EXP3 and EXP4 are shown, where $N^{train}=20, N^{test}=40$. The regions with shaded color around the curves denote its $75\%$ CI.}
        \label{fig:sp_train_test_diffK_40}
\end{figure}

\begin{figure}[H]
\centering
    \includegraphics[width=\columnwidth]{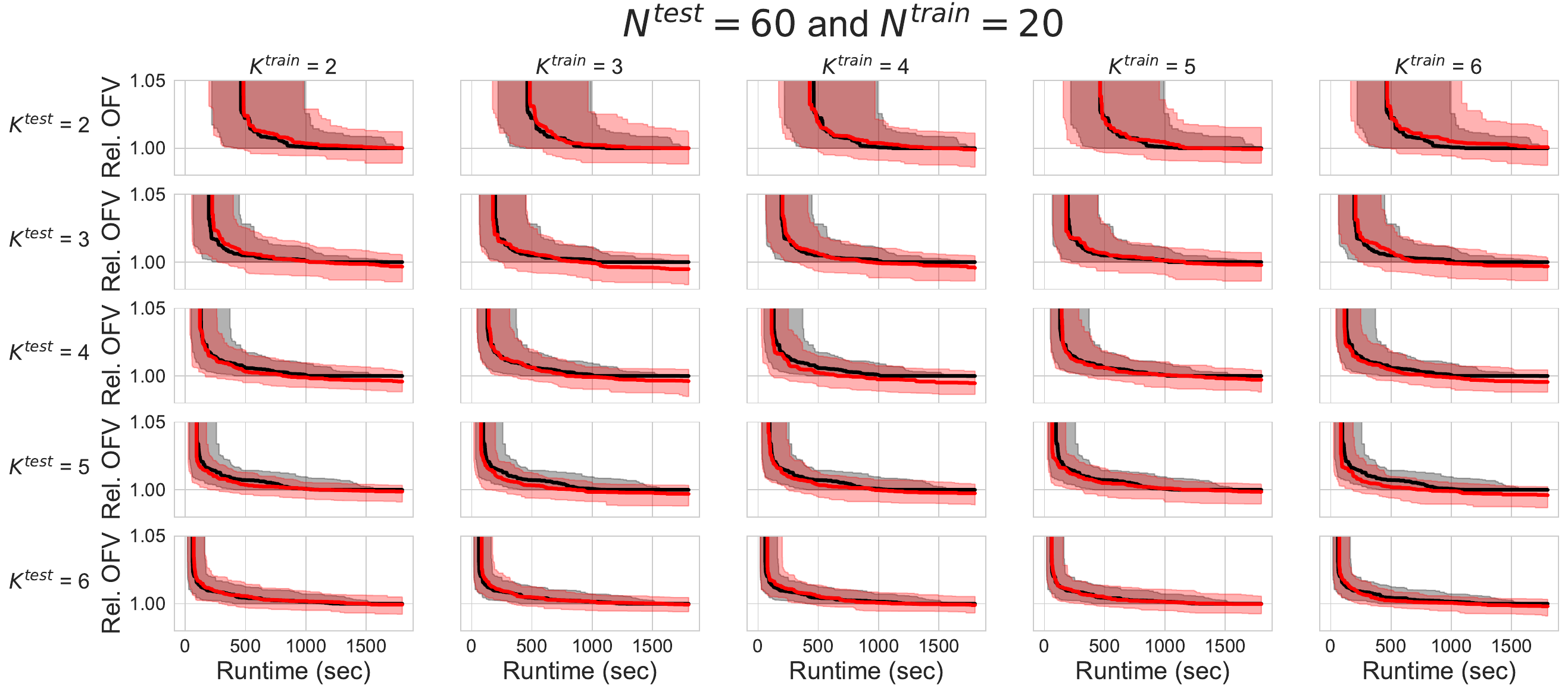}
        \caption{Comparison of results between $K$-B\&B and \textsc{$K$-B\&B-NodeSelection} for 100 instances of the shortest path problem. The results of EXP3 and EXP4 are shown, where $N^{train}=20, N^{test}=60$. The regions with shaded color around the curves denote its $75\%$ CI.}
        \label{fig:sp_train_test_diffK_60}
\end{figure}

\subsection{Mixed problems} \label{app:results_mixed}
\begin{figure}[H]
    \centering
        \includegraphics[width=\columnwidth]{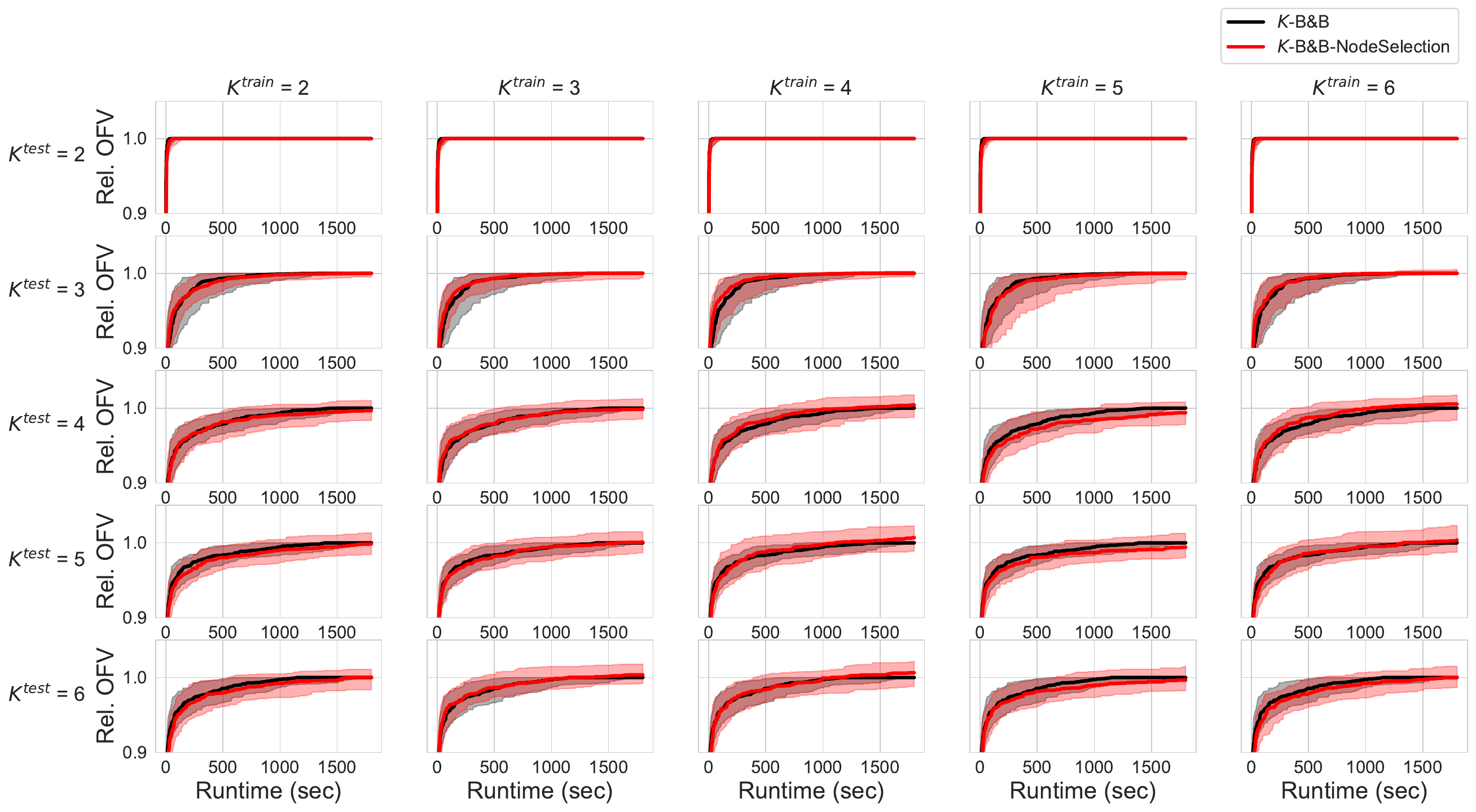}
    \caption{Results of the capital budgeting problem. The ML model that is used is trained on shortest path data.}
    \label{fig:cb_MLSP_train_test_diffK}
\end{figure}

\begin{figure}[H]
    \centering
        \includegraphics[width=\columnwidth]{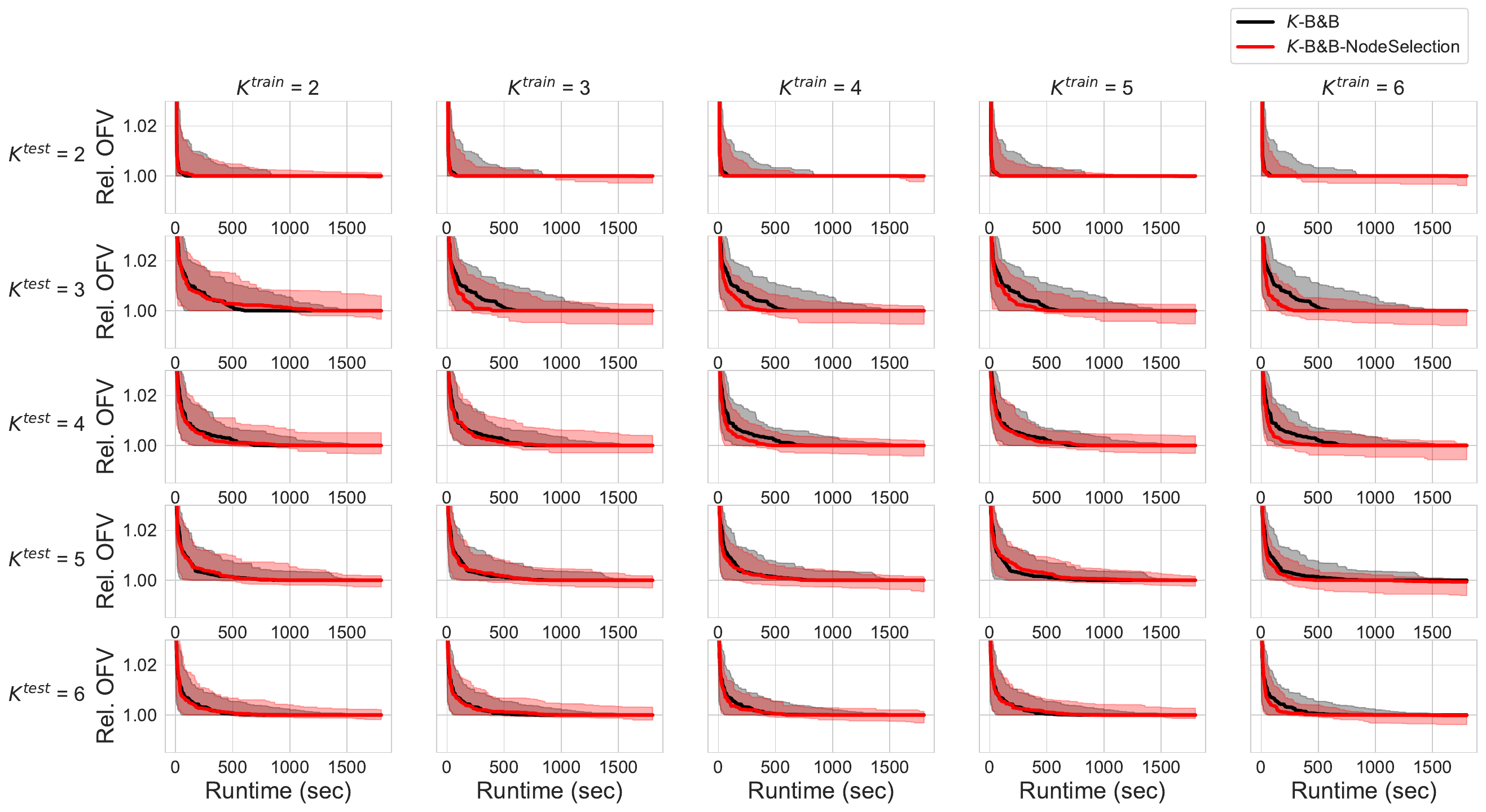}
    \caption{Results of the shortest path problem. The ML model that is used is trained on capital budgeting data.}
    \label{fig:sp_MLCB_train_test_diffK}
\end{figure}

\end{document}